%% file: manuscript.tex
\documentclass[preprint]{elsarticle}

\usepackage{amsmath}
\usepackage{amssymb}
\usepackage{graphicx}
\usepackage{booktabs}
\usepackage[center]{subfigure}
\usepackage{amsmath}
\usepackage{amssymb}
\newcommand{\ds}{\displaystyle}
\def \be{{\bf e}}
\newcommand{\ie}{i.\,e\mbox{.}}
\newcommand{\eg}{e.\,g\mbox{.}}
\newcommand{\cf}{cf\mbox{.}}
\newcommand{\xx}{{\bf x}}
\newcommand{\kk}{{\bf k}}
\newcommand{\Sec}[1]{Section~\ref{#1}}

\newcommand{\Fig}[1]{Fig.~\ref{#1}}

\journal{Computers \& Mathematics with Applications}
\begin{document}
\begin{frontmatter}

\title{Optimization of the Multigrid-Convergence Rate on Semi-structured Meshes by Local Fourier Analysis}

\author[rvt]{B.~Gmeiner\fnref{fn1}}
\ead{bjoern.gmeiner@cs.fau.de}

\author[rvt]{T.~Gradl\fnref{fn1}}
\ead{tobias.gradl@cs.fau.de}

\author[rvt]{F.~Gaspar\fnref{fn2}}
\ead{fjgaspar@unizar.es}

\author[rvt]{U.~R\"ude\fnref{fn1}}
\ead{ulrich.ruede@cs.fau.de}

\address[fn1]{Department of Computer Science, System Simulation, Cauerstr.~11, 91058~Erlangen, Germany}
\address[fn2]{Department of Applied Mathematics, University of Zaragoza. Pedro Cerbuna~12, 50009~Zaragoza, Spain}

%\author{\thanks{Department of Computer Science, System Simulation, Cauerstr.~6, 91058~Erlangen, Germany ({\tt }).} \and
%        Tobias Gradl\thanks{Department of Computer Science, System Simulation, Cauerstr.~6, 91058~Erlangen, Germany} \and
%        Francisco Gaspar\thanks{Department of Applied Mathematics, University of Zaragoza. Pedro Cerbuna~12, 50009~Zaragoza, Spain} \and
%        Ulrich R\"ude\thanks{Department of Computer Science, System Simulation, Cauerstr.~6, 91058~Erlangen, Germany}
%}
%         Paul Duggan\thanks{Composition Department, Society
%         for Industrial and Applied Mathematics, 3600 Univeristy
%         City Science Center, Philadelphia, Pennsylvania,
%         19104-2688 ({\tt duggan@siam.org}).}
%         \and Various A.~U. Thors\thanks{Various Affiliations,
%         supported by various foundation grants.}}

%\baselineskip1cm
\begin{abstract}
In this paper a local Fourier analysis for multigrid methods on tetrahedral grids is presented. Different smoothers for the discretization of the Laplace operator by linear finite elements on such grids are analyzed. A four-color smoother is presented as an efficient choice for regular tetrahedral grids, whereas line and plane relaxations are needed for poorly shaped tetrahedra. A novel partitioning of the Fourier space is proposed to analyze the four-color smoother. Numerical test calculations validate the theoretical predictions. A multigrid method is constructed in a block-wise form, by using different smoothers and different numbers of pre- and post-smoothing steps in each tetrahedron of the coarsest grid of the domain. Some numerical experiments are presented to illustrate the
efficiency of this multigrid algorithm.
\end{abstract}

%42Axx 		Fourier analysis in one variable
%42Bxx 		Fourier analysis in several variables For automorphic theory, see mainly 11F30
%42Cxx 		Nontrigonometric Fourier analysis

%65B99 65F10 65N12 65N22 65N55 

\end{frontmatter}

%\pagestyle{myheadings}

%\markboth{Gmeiner, Gradl, Gaspar, R\"ude}{Optimizing Smoothers on Semi-struct. Meshes by LFA}

\section{\label{sec:introduction}Introduction and terminology}
\input{introduction.tex}

%\section{\label{sec:smoothers}Smoothers for multigrid algorithms}
%\input{smoothers.tex}

%\section{\label{sec:hhg}Hierarchical Hybrid Grids}
%\input{hhg.tex}

\section{\label{sec:lfa}Local Fourier analysis}
\input{lfa.tex}

\section{\label{sec:lfa4color}Four color pattern relaxation}
\input{lfa4color.tex}

\section{\label{sec:LFA_results}Fourier analysis results}
\input{LFA_results.tex}

\section{\label{sec:results}Numerical experiments}
\input{results.tex}

\section{Conclusions}
\input{conclusion.tex}

\section*{Acknowledgments}
This work was partially supported by FEDER/MCYT Projects MTM 2010-16917, the DGA (Grupo consolidado PDIE), 
by the Kompetenznetzwerk f\"ur Technisch-Wissenschaftliches 
Hoch- und H\"ochstleistungsrechnen in Bayern (KONWIHR), and the International
Doctorate Program (IDK) "Identification, Optimization and Control with Applications in Modern Technologies" within the Elite Network of Bavaria.

\input{biblio.tex}

\end{document}

%% file: introduction.tex
Multigrid (MG) methods \cite{Briggs:2000,Hack-1985,TOS01,Wes01} are known to have optimal computational complexity for
solving many numerical problems. However, in practice their performance in solving individual
problems varies significantly. The speed of convergence towards a solution depends on the
numerical properties of the underlying problem, \eg{}, the type of a
differential equation and the method used for discretizing the
equation. Besides that, the user can choose from a variety of algorithms for
the components of the MG method, most prominently the smoother, the
restriction, and the prolongation. Choosing the appropriate components for a
specific problem has also a great impact on the overall performance.

For solving elliptic PDEs, finite element methods are often preferred over
other discretization schemes, because they permit the use of flexible, unstructured
meshes. Algebraic multigrid (AMG) methods \cite{Brandt:1984,Rug87,Stu01,TOS01}
inherently support unstructured meshes by construction. BoomerAMG from the
Hypre package is a popular implementation of AMG \cite{Henson:2002}. Geometric multigrid, in
contrast, relies on a given hierarchy of nested meshes. Geometric multigrid may
achieve a significantly higher performance than algebraic
multigrid in terms of unknowns computed per second.

\begin{figure}
  \caption{\label{fig:2d_levels} Regular refinement example for a
    two-dimensional input mesh. Beginning with the input mesh on the left, each
    successive level of refinement creates a new mesh that has a larger number
    of interior points with structured couplings.}

  \vspace*{1.0ex}
  \centering
  \includegraphics[width=\textwidth]{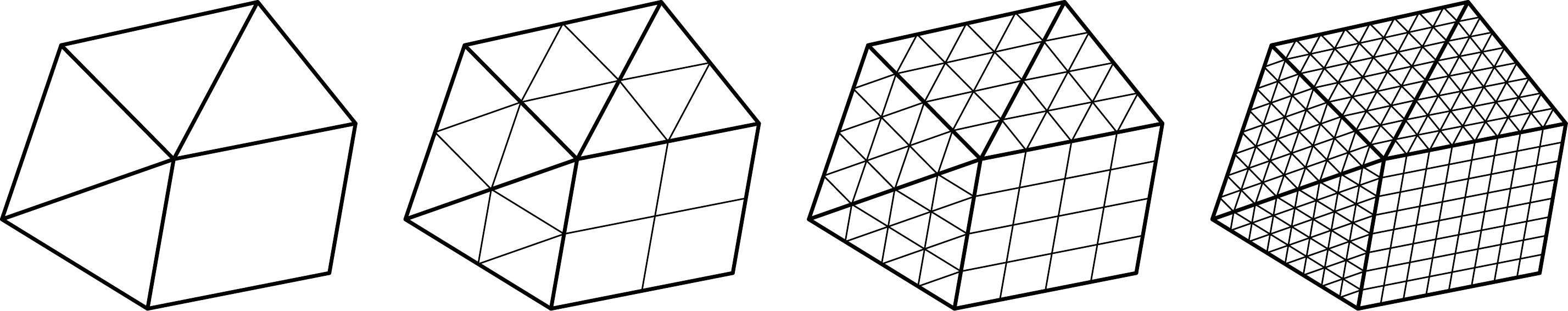}
\end{figure}

The Hierarchical Hybrid Grids (HHG) software framework
\cite{Bergen:2006:CiSE,Gradl} is designed to close this gap between
finite element flexibility and geometric multigrid performance by using a
compromise between structured and unstructured grids. A coarse input finite
element mesh is organized into grid primitives vertices, edges, faces, and
volumes. The primitives are then refined in a structured way, as indicated in
Figure~\ref{fig:2d_levels} for the two-dimensional case. In the case of tetrahedral grids, Bey's refinement strategy
will be considered \cite{Bey95}. There each input tetrahedron is subdivided into eight
child tetrahedra of equal volume, in such a way that each corner of a child coincides to
either a corner or an edge midpoint of the parent. The accuracy of tetrahedral finite elements
w.r.t. their maximal angle is discussed e.g. in \cite{Kri-1992}. The HHG data layout preserves the flexibility of unstructured meshes, while the
regular internal structure of the primitives allows for an efficient
implementation on current computer architectures, especially on parallel
computers. In parallel runs on up to 294\,912 cores, HHG has demonstrated
excellent performance for solving linear systems with up to $10^{12}$
unknowns \cite{Gmeiner:2011:NIC}. Semi-structured meshes also support a local tuning of the
smoothing parameters as will be discussed later in this paper.

Generally the design of MG algorithms and their components should be based on
a careful performance analysis. Such an analysis could be based in the comparison
of convergence rates of different algorithmic variants. For standard cases, such
as the 7-point stencil for the Laplace operator in 3D, convergence rates are
reported in the literature. In our case, however, we will use 15-point stencils
as they arise in the HHG framework with tetrahedral finite elements.

Local Fourier analysis (LFA) \cite{Brandt:1977:MoC} is a very useful tool to predict the asymptotic
convergence factors of MG methods with high accuracy quantitatively. Therefore it is widely used to design efficient
MG algorithms. In the LFA an infinite regular grid is considered and boundary conditions
are ignored. On an infinite grid, the discrete solution and the corresponding error can be represented by
linear combinations of certain complex exponential functions, the Fourier modes, which form a unitary basis of the space of grid functions
with bounded $l_2-$norm. The LFA monograph of Wienands and Joppich \cite{Wienands:2005:PFAfMM}
provides an excellent background for experimenting with Fourier analysis. Recent advances in this context include LFA for hexagonal meshes \cite{Zhou09}, multigrid
as a preconditioner \cite{Wienands:2000:SIAM}, optimal control problems
\cite{Borzi:2007:JCAM}, and discontinuous Galerkin discretizations \cite{Hemker:2004:NLAwA}.
In \cite{Boonen:2008}, an LFA for multigrid methods for the finite element discretization of
the two-dimensional curl-curl equation on a quadrilateral grid has been introduced. 
Recently a generalization of the LFA to triangular grids has been proposed in \cite{Gaspar:2009:SIAM} for a two-dimensional
scalar problem. The key to carrying out this generalization is to express the Fourier transform in a special coordinate system. Our aim in this paper is to show that it is possible to perform an LFA also on 3-dimensional tetrahedral grids so that we can use it to design efficient multigrid solvers in the HHG framework. To choose suitable components of the multigrid method, the LFA will be applied for each
tetrahedron of the input mesh in such a way that the global behavior of the method becomes as efficient as possible.

In a multigrid algorithm, the smoother has the task to reduce the high-frequency
error components, while the coarse-grid correction reduces the low-frequency
components. Both, smoother and coarse-grid correction \cite{de:1996,Moulton-1998} can be tuned to optimally
suit the numerical problem, \ie{} the differential equation and its
discretization. In this paper, we focus on optimizing the smoother depending on
the shape of the tetrahedra.  Besides choosing the
smoother type, two more parameters can be changed to achieve more efficient smoothing
behavior: the under-/over-relaxation parameter ($\omega$) and the number of
pre-/post-smoothing steps per multigrid cycle ($\nu_1$, $\nu_2$). All the
parameters can even be adapted locally to the problem, if the numerical
properties differ strongly across the domain.

\begin{figure}[htb]
\caption{Coloring of a 15-point stencil. Four different colors (and shapes accordingly) are used to break all dependencies between each connected pair of vertices. Planes between the elements are introduced to improve readability. Exemplary two tetrahedra are shown with thicker edges.}
\centering{\includegraphics[height=6cm,width=7cm]{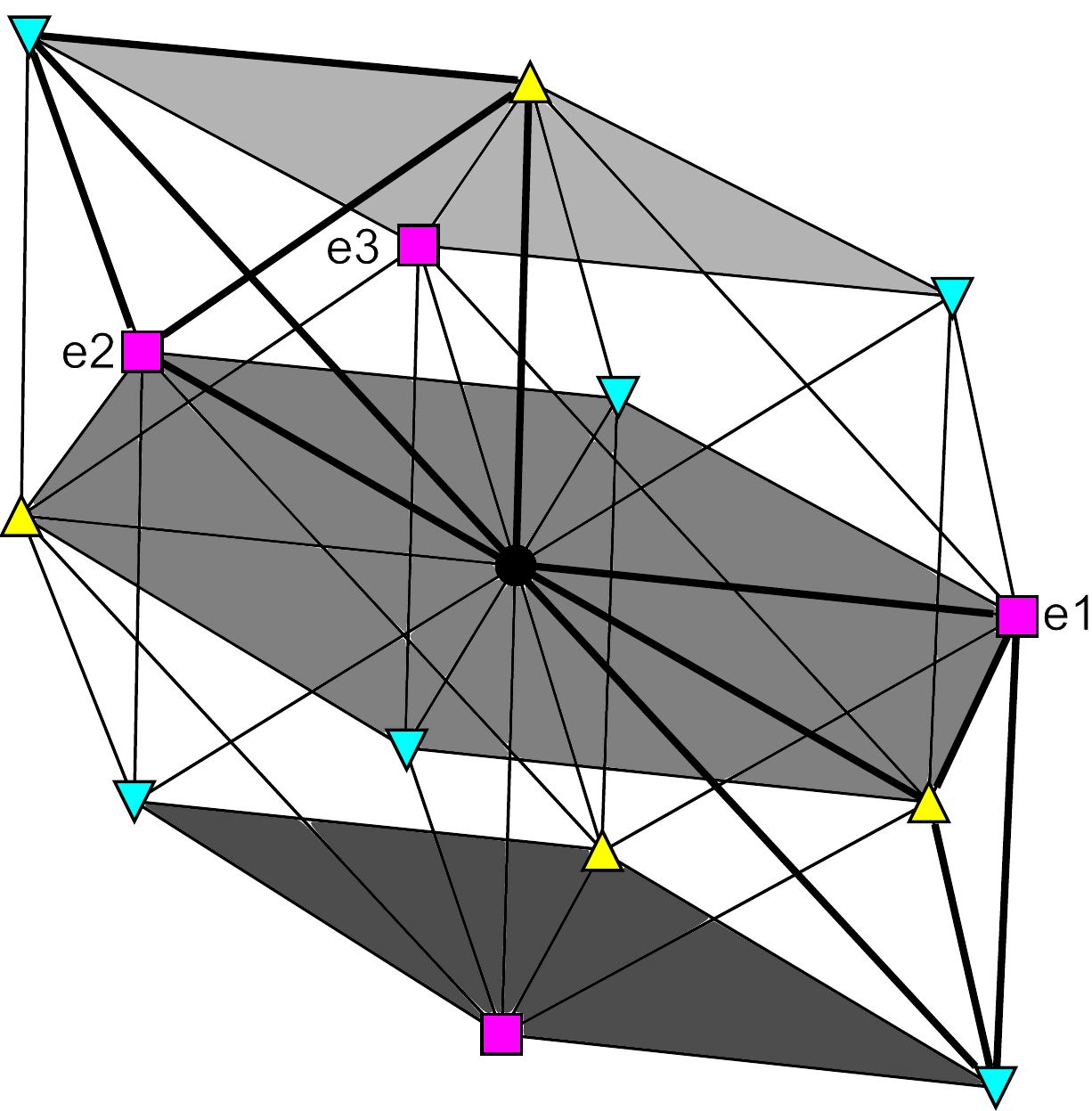}}
\label{four-color}
\end{figure}

The Gauss-Seidel algorithm often is a better smoother than the Jacobi
algorithm, but due to its data dependencies, it is potentially slower.
The data dependencies can be eliminated by coloring the
grid points such that points of the same color are not directly
connected to each other, and can thus be updated in parallel. In a 2D rectangular grid with a 5-point
discretization stencil, e.g. two colors are needed, whereas the 15-point stencil within HHG
solver requires four colors, see Figure~\ref{four-color}. Note that a relaxation parameter $\omega$ can be
chosen individually for each of the colors ($\omega_i = \omega_1, \omega_2, \omega_3, \omega_4$).

In this work, a linear finite element discretization of the Laplace equation will be
analyzed and, as we will see, four-color relaxation results in a very efficient
smoother for regular tetrahedral grids. However, when poorly shaped
tetrahedra occur, then point-wise smoothers are not efficient
anymore. In this case we will consider block-wise smoothers. In particular,
line- and plane-wise smoothers will be used. Note that for Bey
refinement, seven different line smoothers can be defined that correspond with the seven
directions that appear in the connections between each pair of stencil entries,
see Figure~\ref{four-color}. Similarly, seven different plane smoothers can be considered.
Four faces have the orientation of the un-refined tetrahedron face. Three other
face orientations can be spanned up by two vectors connecting two opposing edges. The remaining
multigrid components will be standard, i.e. we will use linear interpolation and its
adjoint as transfer operators and will use Galerkin coarsening for defining the coarse-grid
operators.

The organization of the remaining paper is as follows. In \Sec{sec:lfa} we present the local Fourier analysis for tetrahedral grids. The key here is to consider a basis of the function space appropriate for the structure of the grid.
In \Sec{sec:lfa4color}, a four-color smoother is defined for which smoothing and two-grid analysis will be performed. In this study,
an important aspect is to show the decomposition of the Fourier space into minimal invariant subspaces. These are four-dimensional for the smoothing analysis and have sixteen dimensions for the two-grid analysis.

\Sec{sec:LFA_results} is devoted to Fourier analysis results for the Laplace operator
discretized by linear finite elements on tetrahedral grids. Different smoothers will be
proposed depending on the particular geometry of the grid on each tetrahedral patch. It will be shown that a four-color smoother
will be the best choice for regular tetrahedra. Otherwise block smoothers may be more suitable,
for instance line-smoothers or zebra-plane smoothers.
Finally, in \Sec{sec:results} numerical experiments will show that for domains composed of different types of tetrahedra,
a prescribed global convergence factor, can be obtained using the local LFA  predictions on each tetrahedron
of the input grid.

%% file: lfa.tex
\subsection{General definitions} In this section we extend the local Fourier analysis to multigrid methods on regular tetrahedral grids. The key to
accomplishing this is to use the three-dimensional discrete Fourier transform using coordinates
in non-orthogonal bases fitting the structure of the grid.

\begin{figure}
  \caption{\label{regtet} Regular tetrahedral grid}
  \vspace*{1.0ex}
  \centering
  \includegraphics[scale=0.1]{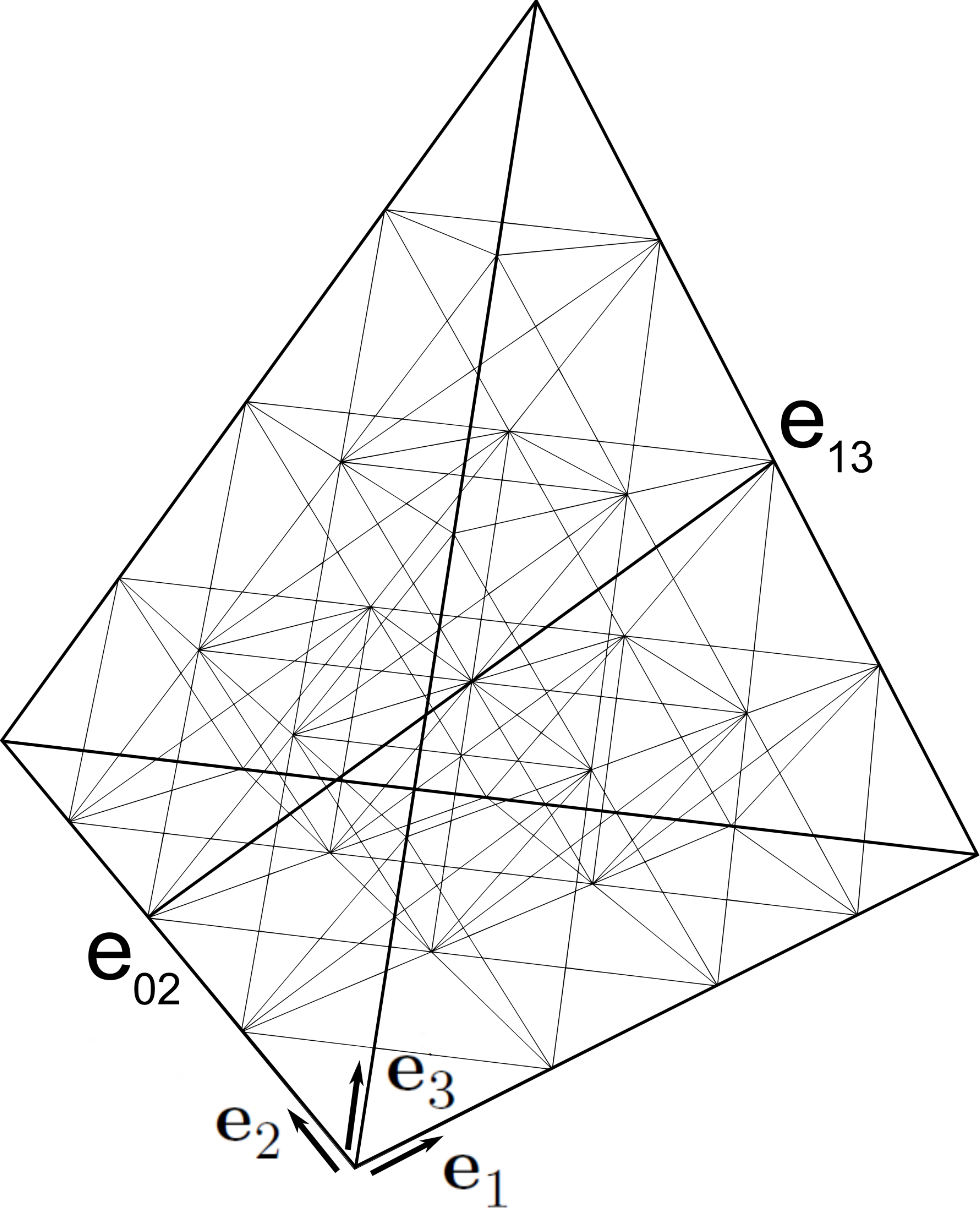}
\end{figure}

Given ${\mathbf h}=(h_1,h_2,h_3)$ a grid spacing, let ${\cal T}_h$ be a regular tetrahedral grid on a fixed coarse tetrahedron ${\cal T}$ and let $\{{\mathbf e}_1, {\mathbf e}_2, {\mathbf e}_3\}$ be unit vectors indicating the direction of three of the
edges of ${\cal T}$; see Figure \ref{regtet}. We extend ${\cal T}_h$ to the infinite grid
\begin{equation}
\label{infinite-grid}
G_h=\{{\mathbf x}= x_1{\mathbf e}_1 + x_2{\mathbf e}_2 + x_3{\mathbf e}_3 | \, x_i = k_i h_i, k_i \in {\bf Z}, i=1,2,3 \},
\end{equation}
such that ${\cal T}_h = G_h \cap {\cal T}$. We write the Fourier frequencies as
${\boldsymbol \theta} = \theta_1 {\mathbf e}_1' + \theta_2 {\mathbf e}_2' +
\theta_3{\mathbf e}_3'$, where $\{{\mathbf e}_1', {\mathbf e}_2', {\mathbf e}_3'  \}$ is the
reciprocal basis to $\{{\mathbf e}_1, {\mathbf e}_2, {\mathbf e}_3\}$, \ie{}, ${\mathbf e}_i \cdot
{\mathbf e}_j' = \delta_{ij}$. In this way, the back Fourier transformation has an expression similar to that in Cartesian coordinates. Therefore, every bounded grid function $u_h({\mathbf x})$ defined on $G_h$ can be written as a formal linear combination of the discrete complex exponential functions called Fourier modes
\begin{equation}
\label{Fourier-mode}
\varphi_h({\boldsymbol \theta}, {\mathbf x}) ={\rm e}^{i (\theta_1 x_1 + \theta_2 x_2 + \theta_3 x_3)},
\end{equation}
where the Fourier frequencies ${\boldsymbol \theta} = \theta_1 {\mathbf e}_1' + \theta_2 {\mathbf e}_2' + \theta_3{\mathbf e}_3'$ vary continuously in ${\bf R}^3$. Since $\varphi_h({\tilde{\boldsymbol \theta}}, {\mathbf x}) = \varphi_h({\boldsymbol \theta}, {\mathbf x})$ for all ${\mathbf x} \in G_h$ if and only if $\tilde{\theta}_i = \theta_i ({\rm mod} \, 2 \pi/h_i), i=1,2,3$, the too high frequencies are not visible on grid $G_h$ due to aliasing effects. These are the Fourier modes with associated frequency ${\boldsymbol \theta}$ such that
for some component $\, |\theta_i| \geq \pi/h_i$. Therefore, the Fourier
space
\[
{\cal F}(G_h) = {\rm span} \{\varphi_h({\boldsymbol \theta},\cdot) \,
\vert \, {\boldsymbol \theta}=(\theta_1,\theta_2,\theta_3) \in {\boldsymbol \Theta}_h\},\]
with ${\boldsymbol \Theta}_h= (-\pi/h_1,\pi/h_1]\times(-\pi/h_2,\pi/h_2]\times(-\pi/h_3,\pi/h_3]$, contains any infinite grid function on $G_h$ with bounded l2-norm. If regular refinement is applied, we
are considering the case of standard multigrid coarsening, $(H_1,H_2,H_3)=(2 h_1,2 h_2, 2 h_3)$, and
therefore the infinite coarse grid $G_H$ is defined as
\[ G_H=\{{\mathbf x}=x_1{\mathbf e}_1 + x_2{\mathbf e}_2 + x_3{\mathbf e}_3 | \, x_i = 2 k_i h_i, k_i \in {\bf Z}, i=1,2,3 \}.
\]
We distinguish high and low frequencies on $G_h$ with respect to $G_{H}$. We
define the subset ${\boldsymbol \Theta}_H$ of low frequencies as $ {\boldsymbol
  \Theta}_H = (-\pi/H_1,\pi/H_1] \times ( -\pi/H_2,\pi/H_2] \times (
    -\pi/H_3,\pi/H_3], $ and the subset of high frequencies ${\boldsymbol
        \Theta}_h \setminus {\boldsymbol \Theta}_H$. This definition is based
      on the fact that only low frequencies are visible on the coarse grid
      $G_H$. Each high-frequency component coincides with a certain
      low-frequency component on $G_H$. In particular we have
      $\varphi_h({\tilde{\boldsymbol \theta}},\cdot) = \varphi_h({\boldsymbol
        \theta},\cdot)$ on $G_H$ with $\tilde{\theta_i} = \theta_i ({\rm
        mod} \, \pi/h_i), i=1,2,3$. \\

The discrete operators considered here are assumed to be linear with constant coefficients.
Neglecting boundary conditions and/or
connections with other neighboring tetrahedra on the coarsest grid, the
discrete problem $L_h u_h = f_h$ can be extended to the whole grid
$G_h$. For a fixed grid point ${\mathbf x}=(x_1,x_2,x_3)\in G_h$, the corresponding equation  reads in stencil notation as
\begin{equation}
L_h u_h (\xx) = \sum_{\kk \in {\cal I}}
s_{\kk}
 u_h(x_1 + k_1 \, h_1, x_2 + k_2 \, h_2, x_3 + k_3 \, h_3)=f_h(\xx),
\label{equ1}
\end{equation}
where $s_{\kk} \in {\mathbf R} $  are constant
coefficients and ${\cal I} \subset {\mathbf Z}^3$ is a finite
index set. From (\ref{equ1}), it can be deduced that the
grid-functions $\varphi_h({\boldsymbol \theta}, \xx)$ given in
(\ref{Fourier-mode}) are formal eigenfunctions of the discrete
operator $L_h$. More precisely, the relation $L_h
\varphi_h({\boldsymbol \theta},\xx) =
\widehat{L_h}({\boldsymbol \theta}) \varphi_h({\boldsymbol
\theta},\xx)$ holds, where
\[
\widehat{L_h}({\boldsymbol \theta})  = \sum_{\kk
\in {\cal I}} s_{\kk} {\rm e }^{i (\theta_1 k_1 h_1+\theta_2 k_2 h_2+\theta_3 k_3 h_3)}
\]
is the corresponding eigenvalue or Fourier symbol of $L_h$.

\subsection{Smoothing analysis}

Recall that one of the basic concepts of multigrid methods is the assumption that high-frequency error components are
smoothed by relaxation, and the low-frequency error components are reduced by coarse grid correction. In Fourier smoothing
analysis, the influence of a smoothing operator $S_h$ on the high-frequency error components is investigated. There an ideal
coarse-grid operator is assumed which annihilates the low frequency error components and leaves the high frequency components unchanged.

Smoothing procedures, such as Gauss-Seidel, line or plane smoothers, have the
property that all $\varphi_h({\boldsymbol \theta},\cdot)$ are eigenfunctions of the smoothing operator. In this case,
the definition of the smoothing factor is given by
\begin{equation}
\mu(S_h) = \sup_{{\boldsymbol \Theta}_h \setminus {\boldsymbol \Theta}_H} \mid \widehat{S_h({\boldsymbol \theta})}\mid
\label{equ:smoothingfactor}
\end{equation}
The situation is more difficult for the four-color relaxation mentioned in \Sec{sec:introduction}, as
the Fourier components are not eigenfunctions of the relaxation operator. \Sec{sec:lfa4color} is devoted to a detailed analysis
of this smoother.

\subsection{Two-grid analysis}

In order to investigate the interplay between relaxation and coarse
grid correction, which is crucial for an efficient multigrid method, it is
necessary to perform a two-grid analysis which takes into account the effect
of transfer operators.  Let $u_h^m$ be an approximation of $u_h$. The error
$e^{m}_h= u_h^m - u_h$ is transformed by a two-grid cycle as $e^{m+1}_h =
M_h^{H} e^m_h$, where $M_h^{H} = S^{\nu_2}_h K_h^{H} S^{\nu_1}_h$ is the
two-grid operator, $K_h^{H} = I_h - P_{H}^h(L_{H})^{-1} R_h^{H} L_h$ the coarse
grid correction operator and $S_h$ is a smoothing operator on $G_h$ with
$\nu_1$ and $\nu_2$ indicating the number of pre- and post-smoothing steps,
respectively. In the definition of $K_h^{H}$, $L_{H}$ is the coarse grid
operator and $P_{H}^h, \, R^{H}_h$ are transfer operators from coarse to fine
grids and vice versa.  The two-grid analysis is the basis for the classical
asymptotic multigrid convergence estimates, and the spectral radius
$\rho(M_h^{H})$ of the operator $M_h^{H}$ indicates the asymptotic convergence
factor of the two-grid method. \\

The Fourier space ${\cal F}(G_h)$ is decomposed into eight-dimensional
subspaces ${\cal F}^8({\boldsymbol \theta}^{000})$, known as $2h$-harmonics, which
are generated by the Fourier modes corresponding to one low-frequency harmonic ${\boldsymbol \theta}^{000} \in {\boldsymbol \Theta}_H$ with coordinates ${\boldsymbol \theta}^{000} = (\theta_1^{000},\theta_2^{000},\theta_3^{000})$ and
seven corresponding high-frequency harmonics in ${\boldsymbol \Theta}_h \setminus {\boldsymbol \Theta}_H$:
\begin{eqnarray}\label{eight-harmonic}
{\cal F}^8({\boldsymbol \theta}^{000}) = {\rm span} \{\varphi_h({\boldsymbol \theta}^{\alpha_1 \alpha_2 \alpha_3},\cdot) \; | \quad \alpha_1,\alpha_2,\alpha_3 \in \{0,1 \}\} \; {\rm with} \; {\boldsymbol \theta}^{000} \in {\boldsymbol \Theta}_H, \nonumber \\
{\boldsymbol \theta}^{\alpha_1 \alpha_2 \alpha_3} = {\boldsymbol \theta}^{000}
-(\alpha_1 {\rm sign}(\theta_1^{000})\pi/h_1,\alpha_2 {\rm sign}(\theta_2^{000})\pi/h_2,\alpha_3 {\rm sign}(\theta_3^{000})\pi/h_3).
\end{eqnarray}
In order to ensure that we deal with non-singular Fourier symbols $\widehat{L_h}({\boldsymbol \theta})$  and $\widehat{L_H}(2 {\boldsymbol \theta})$, we restrict our considerations to the subspace
\[
{\cal F}(G_h) \setminus \bigcup_{{\boldsymbol \theta}^{000} \in {\boldsymbol\Psi}}{\cal F}^8({\boldsymbol \theta}^{000}),
\]
with $\boldsymbol \Psi = \{{\boldsymbol \theta}^{000}
\in \boldsymbol \Theta_H \, \vert \, \widehat{L_H}(2 {\boldsymbol  \theta}^{000})=0 \; {\mbox{\rm
or}} \; \widehat{L_h}({\boldsymbol \theta}^{\alpha_1 \alpha_2 \alpha_3})=0, \ \alpha_1, \alpha_2, \alpha_3 \in\{0,1\} \}.$ The crucial observation is that the coarse-grid correction operator $K_h^{H}$ leaves the spaces of 2h-harmonics invariant for any arbitrary Fourier frequency
${\boldsymbol  \theta}^{000} \in \tilde{\boldsymbol \Theta}_H = {\boldsymbol \Theta}_H \setminus
{\boldsymbol\Psi}$, \ie{}
\[
K_h^H : {\cal
F}^8({\boldsymbol \theta}^{000}) \longrightarrow {\cal
F}^8({\boldsymbol \theta}^{000}).
\]
This same invariance property is
true for most of the smoothers considered here,  \ie{}
$S_h : {\cal F}^8({\boldsymbol \theta}^{000}) \longrightarrow {\cal
F}^8({\boldsymbol \theta}^{000})$. In particular, this is satisfied  for line-type
and plane-type relaxations used in next sections for poorly shaped tetrahedra, see \cite{TOS01, Wienands:2005:PFAfMM} for more details. Therefore, the two-grid operator $M_h^H = S^{\nu_2}_h K_h^H S^{\nu_1}_h$
also leaves the $2h$-harmonic subspaces invariant. Hence, $M_h^H$ is
 equivalent to a block-diagonal matrix consisting of $8
\times 8$ blocks denoted by $\widehat{M_h^H}({\boldsymbol
\theta}^{000}) = M_h^H \mid_{{\cal F}^8({\boldsymbol \theta}^{000})}$,
with ${\boldsymbol \theta}^{000} \in \tilde{{\boldsymbol \Theta}}_H$. For this reason, we can
determine the spectral radius $\rho(M_h^H)$ by calculating the
spectral radius of $(8 \times 8)$-matrices
\begin{equation} \label{factor-coarsening}
\rho(M_h^H) = \sup_{{\boldsymbol \theta}^{000} \in \tilde{{\boldsymbol \Theta}}_H}
\rho(\widehat{M_h^H}({\boldsymbol \theta}^{000})).
\end{equation}
In general, the suprema in~(\ref{equ:smoothingfactor}) and (\ref{factor-coarsening}) are computed numerically by sampling the Fourier space by an enough
fine lattice. All the numerical results presented in this work have been performed by using a grid
with mesh-size $1/32$ in the Fourier space, which provided accurate approximations of the
corresponding factors.

The symbols corresponding to zebra-line and zebra-plane smoothers in tetrahedra grids are
straightforwardly obtained from those in Cartesian coordinates. In fact some of them are identical
and are presented in Appendix A2 in the LFA monograph of Wienands and Joppich~\cite{Wienands:2005:PFAfMM}.

%% file: lfa4color.tex
It is well-known that multicolor relaxations are very efficient smoothers for  multigrid methods.
They are constructed by splitting grid nodes into disjoint sets, with each set having a different color, and
simultaneously updating all nodes of the same color. The best known example of this approach is the red-black
Gauss--Seidel smoother for the five-point Laplace stencil. Such schemes have been
extensively analyzed, see for example \cite{Kuo-1989,Yav95,Yav96}. In this section we focus on the analysis of a four-color relaxation for
discretizations on tetrahedral grids based on 15-point stencils. This smoother is defined in such a way
all unknowns located at grid points with the same color can be updated simultaneously and independently of the
unknowns associated with the other colors, an important property with regards to its parallelization. With this smoother, the two-grid operator
does not leave the $2h$-harmonic subspaces invariant, and therefore a more elaborate Fourier analysis is required. \\

\subsection{Smoothing analysis}
In this section we focus in the analysis of the four-color relaxation. The infinite grid $G_h$ (\ref{infinite-grid}) is subdivided into four different type of grids $G_h^i, \; i=0,1,2,3$, each associated with a different color,
$$
G^i_h=\{\xx=(x_1,x_2,x_3)
\, \vert
\, x_j=k_jh_j, \, \, k_j\in {\bf Z}, \, j=1,2,3 , \ k_1+k_2+k_3=i \, ({\rm mod \, } 4) \},
$$
where $(x_1,x_2,x_3)$ are the coordinates of the nodes of the grid in the basis $\{\be_1,\be_2,\be_3\}$.
The special choice of this basis gives a different color
to each vertex of  a tetrahedron,  in such way that the unknowns with the same color
have no direct connection with each other, see Figure \ref{four-color}. More concretely,
for any grid point $\xx=(x_1,x_2,x_3) \in G_h^i$ only  the grid
points of the form $(x_1 + \kappa_1 h_1, x_2 + \kappa_2 h_2, x_3 + \kappa_3 h_3)$, with
\begin{eqnarray*}
(\kappa_1,\kappa_2,\kappa_3) &\in& \{(1,0,0),(1,1,0),(0,1,0),(-1,0,0),(-1,-1,0),(0,-1,0), \\
 & &(-1,0,-1), (0,0,-1),(-1,-1,-1),(0,-1,-1), \\ && (0,0,1),(1,0,1),(1,1,1),(0,1,1)\}
\end{eqnarray*}
contribute to the corresponding stencil (apart from $\xx$). Note that none of these points belong to the grid $G_h^i$, i.e., $\kappa_1 + \kappa_2 + \kappa_3 \neq 0 \, ({\rm mod \, } 4)$. The four--color smoother is defined as a sequence of partial Jacobi sweeps in
which the variables are updated, for example, in the order $G^0_h, G^1_h,G^2_h$ and
$G^3_h$. Notice that in each partial relaxation step, only the grid points of $G^i_h$ are processed, whereas
the remaining points are not treated. The smoothing operator $S_h^i(\omega_i)$ associated with the $i^{th}$--partial step is given by
\begin{equation}
S^i_h(\omega_i) v_h(\xx)=\left \{ \begin{array}{cl} [(I_h-\omega_i D_h^{-1}L_h)
v_h](\xx), & \xx \in G_h^i, \\
v_h(\xx), & \xx \in G_h \backslash G_h^i,
\end{array}
\right.
\label{operatorS}
\end{equation}
where $D_h$ is the diagonal part of the discrete operator $L_h$ and $\omega_i$ is a damping parameter related to the corresponding color.
Therefore the complete four-color smoothing operator is given by the product of these four partial step operators,
$S_h=S^3_h(\omega_3) S^2_h(\omega_2) S^1_h(\omega_1) S^0_h(\omega_0)$.

\begin{figure}
  \caption{\label{fig:3dfrequencies} Four-dimensional invariant harmonic subspaces for the four-color smoother}

  \vspace*{1.0ex}
  \centering
  \includegraphics[width=\textwidth]{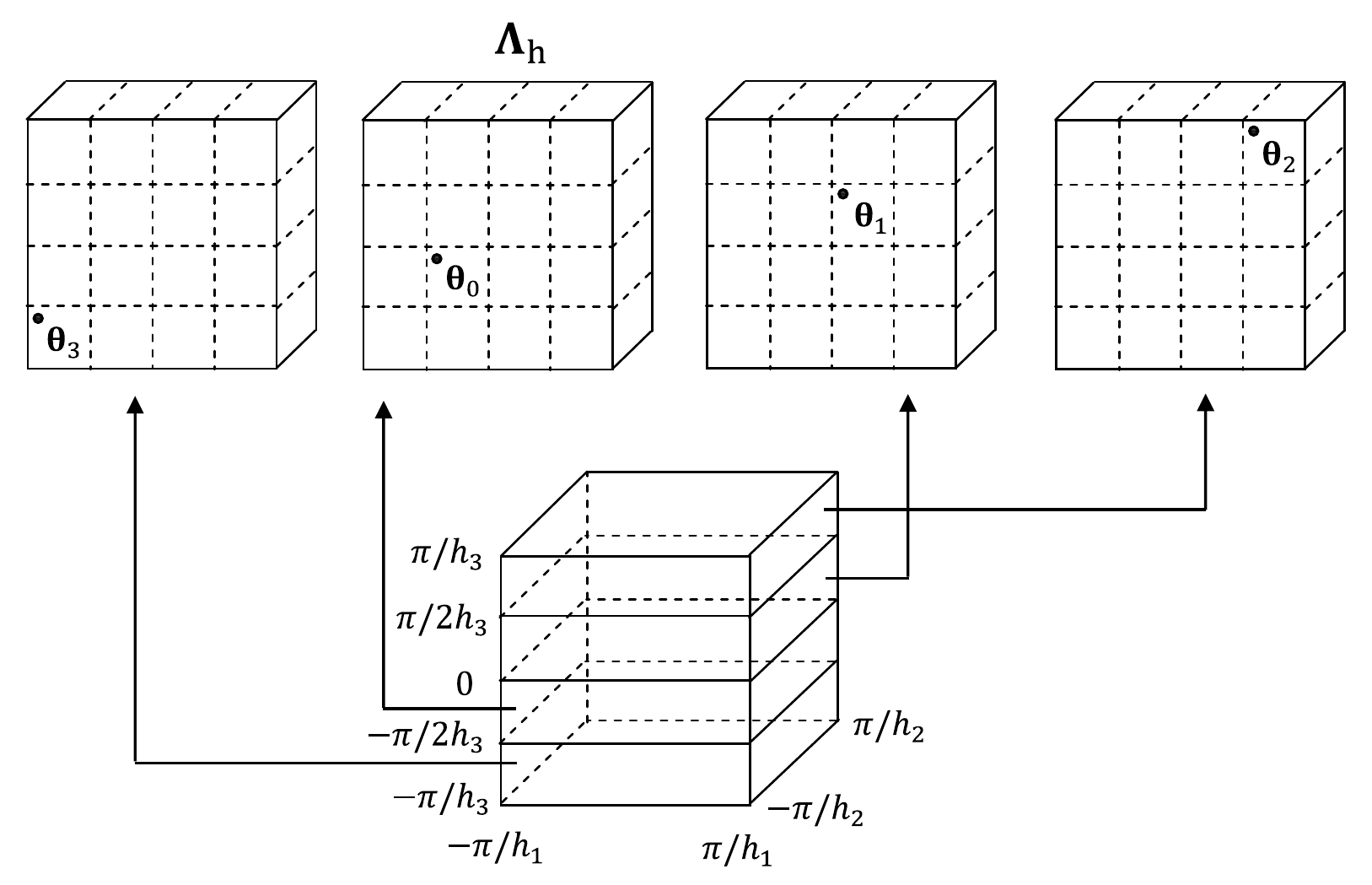}
\end{figure}

%\subsection{Smoothing analysis}
As the Fourier modes  $\varphi_h({\boldsymbol \theta}, \cdot)$ are not eigenfunctions of the smoothing operators $S_h^i(\omega_i), i=0,\ldots,3$, a new strategy is necessary for the computation of the smoothing factor. For this, let us define the set of frequencies $\boldsymbol \Lambda_h$ given by
$$ \boldsymbol \Lambda_h=\left(-\ds\frac{\pi}{h_1},\ds\frac{\pi}{h_1}\right]\times\left(-\ds\frac{\pi}{h_2},\ds\frac{\pi}{h_2}\right]
\times\left(-\ds\frac{\pi}{2h_3},0\right].
$$
Given a frequency ${\boldsymbol \theta}_0 \in \boldsymbol \Lambda_h$, let us consider the four-dimensional subspace
${\cal F}^4(\boldsymbol \theta_{0})$ given by
\begin{equation}\label{3ar} {\cal F}^4(\boldsymbol \theta_{0})={\rm span}\{\varphi_h(\boldsymbol
\theta_{0}, \cdot),\varphi_h(\boldsymbol
\theta_{1}, \cdot),\varphi_h(\boldsymbol
\theta_{2}, \cdot),\varphi_h(\boldsymbol
\theta_{3}, \cdot)\},\end{equation}
where
$\boldsymbol \theta_{\sigma}=(\boldsymbol \theta_{0}
+\frac{\pi\sigma}{2}{\bf h}^{-1}) ({\rm mod \,} 2\pi{\bf h}^{-1}), \ \sigma=1,2,3$ and ${\bf h}^{-1}=(h_1^{-1},h_2^{-1},h_3^{-1})$. From the definition of $\boldsymbol \theta_{\sigma}$, it results that if $\boldsymbol \theta_{0} \in \boldsymbol \Lambda_h$, each of the four coupled frequencies belongs to a different subregion of the frequency domain ${\boldsymbol \Theta}_h$, as observed in Figure \ref{fig:3dfrequencies}. If $\boldsymbol \theta_{0}$ covers $\boldsymbol \Lambda_h$, the four frequencies defining the Fourier modes generating ${\cal F}^4(\boldsymbol \theta_{0})$ fill the whole
frequency domain ${\boldsymbol \Theta}_h$. Due to this fact, the Fourier space is given by the direct sum of subspaces ${\cal F}^4(\boldsymbol \theta_{0})$, i.e.,
$${\cal F}(G_h) = \displaystyle \bigoplus_{\boldsymbol \theta_{0} \in \boldsymbol \Lambda_h} {\cal F}^4(\boldsymbol \theta_{0}).$$
We now show that for each ${\boldsymbol \theta}_0 \in \boldsymbol \Lambda_h$, the operators $S^i_h(\omega_i), i=0,\ldots,3$ leave the spaces ${\cal F}^4(\boldsymbol \theta_{0})$ invariant, i.e.,
$$
S_h^i(\omega_i): {\cal F}^4(\boldsymbol \theta_{0}) \to {\cal F}^4(\boldsymbol \theta_{0}), \quad i=0,\ldots,3, \quad \boldsymbol \theta_{0}\in \boldsymbol \Lambda_h.
$$
If operator $S_h^0(\omega_0)$ is applied to $\varphi_h(\boldsymbol \theta_{0},\xx)$, we have that
\begin{equation}
S^0_h(\omega_0) \varphi_h(\boldsymbol \theta_{0},\xx)=\left \{ \begin{array}{cl} \lambda_0
\varphi_h(\boldsymbol \theta_{0},\xx), & \xx \in G_h^0, \\
\varphi_h(\boldsymbol \theta_{0},\xx), & \xx \in G_h \backslash G_h^0,
\end{array}
\right.
\label{operatorSeigen}
\end{equation}
where $\lambda_{0}=1-\omega_0 \widehat{D_h}^{-1} \widehat{L_h} (\boldsymbol
\theta_{0})$ is the Fourier symbol of Jacobi relaxation, see \cite{TOS01}, with $\widehat{D_h}$ and $\widehat{L_h}(\boldsymbol
\theta_{0})$  the symbols of the operators $D_h$ and $L_h$ respectively.
Expression (\ref{operatorSeigen}) is not yet a Fourier representation of $S_h^0$, however
it is easy to see that $S^0_h(\omega_0) \varphi_h(\boldsymbol
\theta_{0},\xx) \in {\cal F}^4(\boldsymbol \theta_{0})$, i.e.,
\begin{equation} \label{coef-4arm}
\begin{array}{l}
S^0_h(\omega_0) \varphi_h(\boldsymbol
\theta_{0},\xx)=\alpha_{0}(\boldsymbol
\theta_{0}) \varphi_h(\boldsymbol
\theta_{0},\xx)+\alpha_{1} (\boldsymbol \theta_{0})
\varphi_h(\boldsymbol \theta_{1},\xx)+\alpha_{2} (\boldsymbol
\theta_{0})
\varphi_h(\boldsymbol \theta_{2},\xx) + \\ \quad \alpha_{3} (\boldsymbol
\theta_{0})
\varphi_h(\boldsymbol \theta_{3},\xx)
\end{array}
\end{equation}
with
$$
\alpha_{0}(\boldsymbol
\theta_{0})=\ds\frac{\lambda_{0}+3}{4}, \quad \alpha_{1}(\boldsymbol
\theta_{0})=\alpha_{2}(\boldsymbol
\theta_{0})=\alpha_{3}(\boldsymbol
\theta_{0})=\ds\frac{\lambda_{0}-1}{4}.
$$
These coefficients are obtained by evaluating expressions (\ref{operatorSeigen}) and (\ref{coef-4arm}) in the
four subgrids $G_h^i, i=0,\ldots,3$:
$$
\begin{array}{rl}
\alpha_{0}(\boldsymbol
\theta_{0})+\alpha_{1}(\boldsymbol
\theta_{0})+\alpha_{2}(\boldsymbol
\theta_{0})+\alpha_{3}(\boldsymbol
\theta_{0})=\lambda_{0}, & \hbox{if } \xx \in G_h^0, \\
\alpha_{0}(\boldsymbol
\theta_{0})+e^{i \pi /2} a_{1}(\boldsymbol
\theta_{0})+e^{i \pi } a_{2}(\boldsymbol
\theta_{0})+ e^{i 3 \pi /2} a_{3}(\boldsymbol
\theta_{0})=1, & \hbox{if }\xx \in G_h^1, \\
\alpha_{0}(\boldsymbol
\theta_{0})+e^{i \pi} a_{1}(\boldsymbol
\theta_{0})+e^{i 2 \pi } a_{2}(\boldsymbol
\theta_{0})+ e^{i 3 \pi} a_{3}(\boldsymbol
\theta_{0})=1, & \hbox{if }\xx \in G_h^2, \\
\alpha_{0}(\boldsymbol
\theta_{0})+e^{i 3 \pi /2} a_{1}(\boldsymbol
\theta_{0})+e^{i 3 \pi } a_{2}(\boldsymbol
\theta_{0})+ e^{i 9 \pi /2} a_{3}(\boldsymbol
\theta_{0})=1, & \hbox{if }\xx \in G_h^3,
\end{array}
$$
and solving the corresponding $4\times4$-system. Analogously, if we denote by
\begin{equation}
S^0_h(\omega_0) \varphi_h(\boldsymbol \theta_{\sigma},\xx)=\left \{
\begin{array}{cl} \lambda_{\sigma} \varphi_h(\boldsymbol
\theta_{\sigma},\xx)
, & \xx \in G_h^0, \\
\varphi_h(\boldsymbol \theta_{\sigma},\xx), & \xx \in G_h
\backslash G_h^0,
\end{array}
\right.
\end{equation}
where $\lambda_{\sigma}=1-\omega_{\sigma} \widehat{D_h}^{-1} \widehat{L_h} (\boldsymbol
\theta_{\sigma})$ for $\sigma=1,2,3$, we can write
\begin{equation} \label{coef-4arm2}
\begin{array}{l}
S^0_h(\omega_0) \varphi_h(\boldsymbol
\theta_{\sigma},\xx)=\alpha_{0}(\boldsymbol
\theta_{\sigma}) \varphi_h(\boldsymbol
\theta_{0},\xx)+\alpha_{1} (\boldsymbol \theta_{\sigma})
\varphi_h(\boldsymbol \theta_{1},\xx)+\alpha_{2} (\boldsymbol
\theta_{\sigma})
\varphi_h(\boldsymbol \theta_{2},\xx) + \\ \quad \alpha_{3} (\boldsymbol
\theta_{\sigma})
\varphi_h(\boldsymbol \theta_{3},\xx), \quad \sigma=1,2,3,
\end{array}
\end{equation}
with
$$
\alpha_{\sigma'}(\boldsymbol
\theta_{\sigma})=\left\{\begin{array}{ll}
\ds\frac{\lambda_{\sigma}+3}{4}, & \hbox{if } \sigma'=\sigma,\\
\ds\frac{\lambda_{\sigma}-1}{4}, & \hbox{otherwise.}
\end{array}
\right.
$$
We conclude that for each ${\boldsymbol \theta}_0 \in \boldsymbol \Lambda_h$, the operator $S^0_h(\omega_0)$ leaves the space ${\cal F}^4(\boldsymbol \theta_{0})$ invariant, and its matrix representation is the $(4\times 4)$--matrix
\begin{equation} \label{sh0}
\widehat{S^0_h}(\boldsymbol \theta_{0})=\frac{1}{4} \left
(\begin{array}{cccc} \lambda_0+3 &\lambda_1-1 & \lambda_2-1 & \lambda_3-1 \\
 \lambda_0-1 & \lambda_1+3 & \lambda_2-1 & \lambda_3-1 \\
 \lambda_0-1 &\lambda_1-1 & \lambda_2+3 & \lambda_3-1 \\
 \lambda_0-1 &\lambda_1-1 & \lambda_2-1 & \lambda_3+3
\end{array} \right).
\end{equation}
Analogously, the $(4\times 4)$--matrices $\widehat{S^1_h}(\boldsymbol
\theta_{0})$, $\widehat{S^2_h}(\boldsymbol \theta_{0}), \widehat{S^3_h}(\boldsymbol \theta_{0})$ associated with
the smoothing operators $S^1_h(\omega_1)$, $S^2_h(\omega_2)$ and $S^3_h(\omega_3)$ can be determined, yielding to
the following expressions:
\begin{eqnarray}
\widehat{S^1_h}(\boldsymbol \theta_{0})&=&\frac{1}{4} \left
(\begin{array}{cccc} \lambda_0+3 &i(\lambda_1-1) & 1-\lambda_2 & i(1-\lambda_3) \\
 i(1-\lambda_0) & \lambda_1+3 & i(\lambda_2-1) & 1-\lambda_3 \\
 1 - \lambda_0 &i(1-\lambda_1) & \lambda_2+3 & i(\lambda_3-1) \\
 i(\lambda_0-1) & 1-\lambda_1 & i(1-\lambda_2) & \lambda_3+3
\end{array} \right), \label{sh1}\\
\widehat{S^2_h}(\boldsymbol \theta_{0})&=&\frac{1}{4} \left
(\begin{array}{cccc} \lambda_0+3 &1-\lambda_1 & \lambda_2-1 & 1-\lambda_3 \\
 1-\lambda_0 & \lambda_1+3 & 1-\lambda_2 & \lambda_3-1 \\
 \lambda_0-1 &1-\lambda_1 & \lambda_2+3 & 1-\lambda_3 \\
 1-\lambda_0 &\lambda_1-1 & 1-\lambda_2 & \lambda_3+3
\end{array} \right), \label{sh2} \\
\widehat{S^3_h}(\boldsymbol \theta_{0})&=&\frac{1}{4} \left
(\begin{array}{cccc} \lambda_0+3 &i(1-\lambda_1) & 1-\lambda_2 & i(\lambda_3-1) \\
 i(\lambda_0-1) & \lambda_1+3 & i(1-\lambda_2) & 1-\lambda_3 \\
 1 - \lambda_0 &i(\lambda_1-1) & \lambda_2+3 & i(1-\lambda_3) \\
 i(1-\lambda_0) & 1-\lambda_1 & i(\lambda_2-1) & \lambda_3+3
\end{array} \right). \label{sh3}
\end{eqnarray}
Finally, we have that $\widehat{S_h}(\boldsymbol \theta_{0})=\widehat{S_h^3}(\boldsymbol \theta_{0})\widehat{S_h^2}(\boldsymbol \theta_{0})\widehat{S_h^1}(\boldsymbol \theta_{0})\widehat{S_h^0}(\boldsymbol \theta_{0})$ is the $(4\times 4)$--matrix
representation of $S_h$ in ${\cal F}^4(\boldsymbol \theta_{0})$ and the
smoothing factor per sweep, for $\nu$ consecutive sweeps, is given
by
\begin{equation} \label{factor-smoothing2}
\mu(S_h,\nu)=\sup_{\boldsymbol{\theta}_{0}\in \boldsymbol \Lambda_h} \,
\sqrt[\nu]{\rho (\widehat{Q}(\boldsymbol \theta_{0})
\widehat{S_h}^\nu(\boldsymbol \theta_{0}))},
\end{equation}
where
 $\widehat{Q}(\boldsymbol \theta_{0})=
 {\rm diag}(q(\boldsymbol{\theta}_{0}),
 q(\boldsymbol{\theta}_{1}),
 q(\boldsymbol{\theta}_{2}),
 q(\boldsymbol{\theta}_{3}))$ is a diagonal matrix with
 $$
\begin{array}{l}
q(\boldsymbol{\theta}_{\sigma})=\left\{\begin{array}{ll}0, & \hbox{if }
\boldsymbol{\theta}_{\sigma} \in\boldsymbol \Theta_H, \\ 1, & \hbox{otherwise,}
\end{array}\right.
\end{array}
$$
for $\sigma=0,\ldots,3$. Notice that the projection matrix $\widehat{Q}(\boldsymbol \theta_{0})$ annihilates the low--frequency error components and leaves the
high--frequency error components unchanged.

\subsection{Two-grid analysis}
The smoothing factor give by formula \ref{factor-smoothing2} should accurately approximate the actual multigrid convergence, if the considered multigrid coarse-grid correction is close to the
idealized coarse-grid correction. In other case, a two-grid analysis is necessary to
understand the interplay between the smoothing and the coarse-grid correction process.
To perform a two-grid analysis, Fourier space ${\cal F}(G_h)$ is divided into a direct sum of sixteen--dimensional subspaces
${\cal F}(G_h) = \oplus {\cal F}^{16}(\boldsymbol \theta_{0}^{000})$, where $\boldsymbol \theta_{0}^{000} \in \boldsymbol \Lambda_H$, with
$$ \boldsymbol \Lambda_H=\left(-\ds\frac{\pi}{2h_1},\ds\frac{\pi}{2h_1}\right]\times\left(-\ds\frac{\pi}{2h_2},\ds\frac{\pi}{2h_2}\right]
\times\left(-\ds\frac{\pi}{2h_3},0\right].
$$
These sixteen-dimensional subspaces are defined as
\[
{\cal F}^{16} ({\boldsymbol \theta_{0}^{000}})={\cal
F}^{4}({\boldsymbol \theta_{0}^{000}})\oplus {\cal
F}^{4}({\boldsymbol \theta_{0}^{110}})\oplus {\cal
F}^{4}({\boldsymbol \theta_{0}^{100}}) \oplus {\cal
F}^{4}({\boldsymbol \theta_{0}^{010}}),
\]
where
${\boldsymbol \theta_{0}^{110}}, {\boldsymbol \theta_{0}^{100}}, {\boldsymbol
\theta_{0}^{010}}$ are the three frequencies associated with
${\boldsymbol \theta_{0}^{000}}$ by means of (\ref{eight-harmonic}) and the subspace ${\cal
F}^4({\boldsymbol \theta_0^{ij0}})$ is generated as in (\ref{3ar}).

\begin{figure}
  \caption{\label{fig:3dfrequencies2} Sixteen-dimensional minimal invariant harmonic subspaces for the four-color smoother}

  \vspace*{1.0ex}
  \centering
  \includegraphics[width=\textwidth]{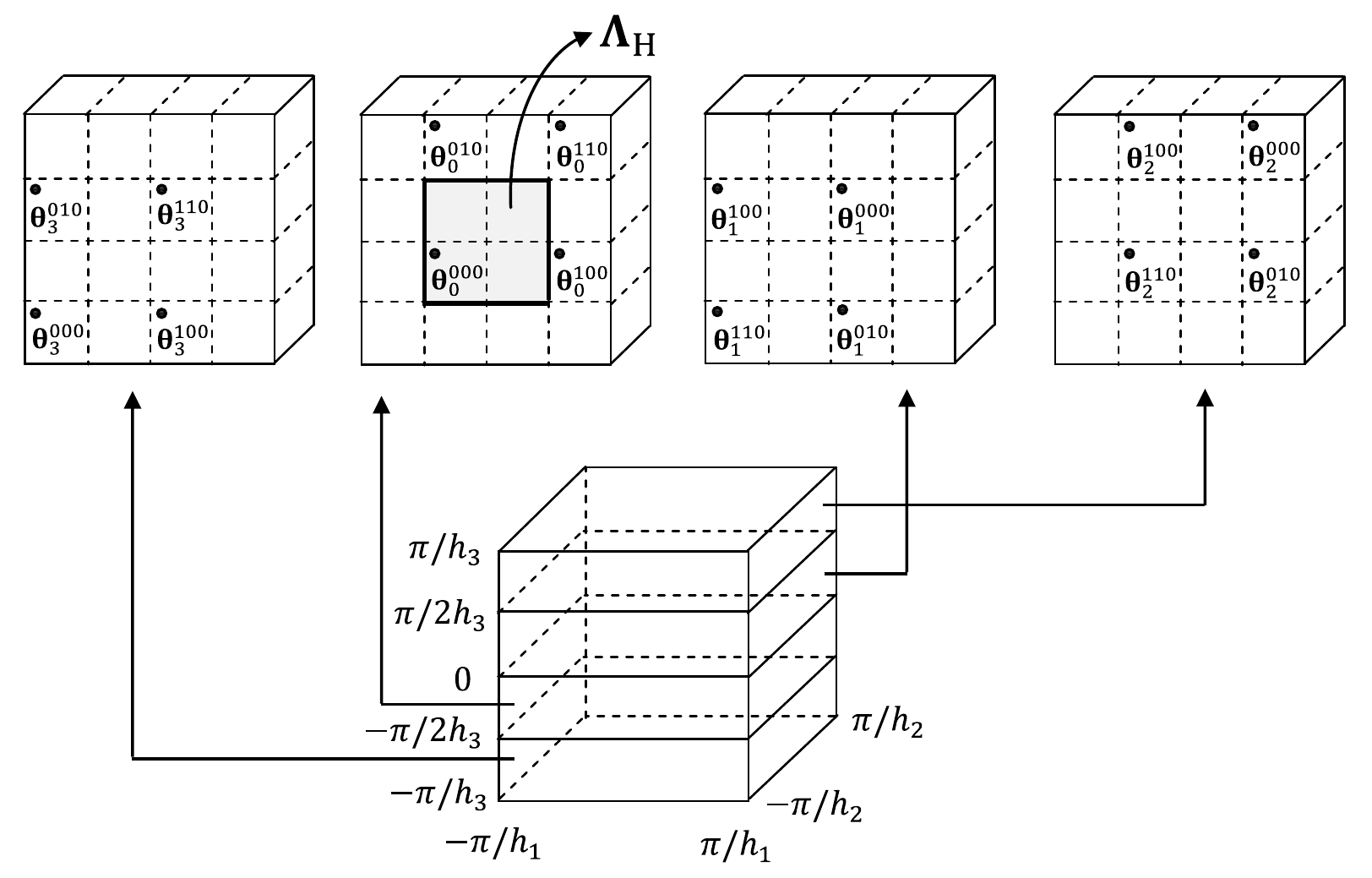}
\end{figure}

 Figure \ref{fig:3dfrequencies2} illustrates the location of the sixteen different frequencies $\boldsymbol \theta^{ij0}_{\sigma}$ and the region $\boldsymbol \Lambda_H$. From the invariance results of the previous sections, it can be observed that
the subspaces $ {\cal F}^{16} ({\boldsymbol \theta_{0}^{000}})$ are invariant for both the four--color
smoother and the standard coarsening operator. Notice that in each subspace $ {\cal F}^{16} ({\boldsymbol \theta_{0}^{000}})$, we can find
associated two low-frequencies ${\boldsymbol \theta_{0}^{000}}$ and ${\boldsymbol \theta_{1}^{000}}$, and their seven corresponding high-frequency harmonics. Therefore, if ${\boldsymbol \theta_{0}^{000}}$
covers the whole $\Lambda_H$ then two coupled low-frequencies fill the subset ${\boldsymbol \Theta}_H$
of low frequencies. This fact together with the coupling of a low-frequency with seven high frequency modes implies that the sixteen frequencies defining the Fourier modes generating $ {\cal F}^{16} ({\boldsymbol \theta_{0}^{000}})$ fill the whole frequency domain ${\boldsymbol \Theta}_h$.

The matrix representation of $S^i_h(\omega_i)$  in ${\cal F}^{16}(\boldsymbol \theta_{0}^{000}), \, \boldsymbol \theta_{0}^{000} \in \boldsymbol \Lambda_H$, is now the $(16\times 16)$--block matrix ${\rm diag}[\widehat{S^i_h}(\boldsymbol \theta_{0}^{000}),\widehat{S^i_h}(\boldsymbol \theta_{0}^{110}),
\widehat{S^i_h}(\boldsymbol \theta_{0}^{100}),\widehat{S^i_h}(\boldsymbol \theta_{0}^{010})]$, where the $4\times4$ diagonal matrices  $\widehat{S^i_h}(\boldsymbol \theta_{0}^{ij0})$ are given in (\ref{sh0})-(\ref{sh3}). The matrix representation $\widehat{K^H_h}(\boldsymbol \theta_{0}^{000})$  of the coarse-grid correction operator $K^H_h$ in ${\cal F}^{16}(\boldsymbol \theta_{0}^{000})$ is a $(16\times 16)$--block matrix, constructed in the usual way from the two low-frequencies $\theta_0^{000}$ and $\theta_1^{000}$ defining the two smooth Fourier modes in $F^{16}(\theta_0^{000})$.

The definition of the asymptotic convergence factor for a two-grid cycle with $\nu_1$ pre--smoothing and
$\nu_2$ post--smoothing steps now becomes
\begin{equation} \label{factor-coarsening2}
\rho(M^H_h) =\sup_{\boldsymbol{\theta_{0}^{000}}\in \tilde{\boldsymbol \Lambda}_{H}} \,
\rho ( \widehat{M^H_h}(\boldsymbol \theta_{0}^{000}))=\sup_{\boldsymbol{\theta_{0}^{000}}\in \tilde{\boldsymbol \Lambda}_{H}} \,
\rho ( \widehat{S_h}^{\nu_2}(\boldsymbol \theta_{0}^{000})
\widehat{K^H_h} (\boldsymbol \theta_{0}^{000})
\widehat{S_h}^{\nu_1}(\boldsymbol \theta_{0}^{000}) ),
\end{equation}
with
$$
\begin{array}{r}
\tilde{\boldsymbol \Lambda}_{H}=\boldsymbol \Lambda_{H} \backslash \{{\boldsymbol \theta}_{0}^{000}
\in \boldsymbol \Lambda_{H} \, \vert \, \widehat{L_H}(2 {\boldsymbol  \theta}^{jkl}_{\sigma})=0,  \; {\mbox{\rm
or}}\\ \; \widehat{L_h}({\boldsymbol \theta}^{jkl}_{\sigma})=0,  \ \sigma=0,1,2,3, \; j,k,l\in\{0,1\}\}.
\end{array}
$$ 

%% file: LFA_results.tex
This section is devoted to show some LFA results in order to design efficient
multigrid methods for different shapes of tetrahedra. We mainly focus on
the choice of good smoothers for poorly shaped tetrahedra which are inevitably produced
by commonly used mesh generators.

In order to compare the convergence factors predicted by LFA, a W-cycle is used
since a two-grid cycle is often not practical since it would require a high computational
effort on the coarse grid. However,the convergence factors are  similar to a two-grid cycle.
The LFA assumes an infinitely large regular domain. In practice, such a domain
can be approximated by a regular grid with a very large number of unknowns. In
experiments, structured meshes of $2.1 \cdot 10^{6}$ elements (129 grid points in each direction) easily
turned out to be sufficient to exhibit the asymptotic behaviour. Similar number of grid points in each directions are necessary in 2D. 
In the experiment, the solution is initialized with random values between $-1.0$ and $1.0$. The errors
for determining the convergence rates are measured with the discrete L2 norm, using the observed rate of error reduction
in the 151st cycle as approximation to the asymptotic value.

If not stated differently, all relaxation parameters are found by LFA for the remaining paper.
At first we are interested to find optimal relaxation parameters 
w.r.t. the convergence rate. To cover the complete search space, all different 
relaxation parameters up to $\omega=2.0$ are predicted with a resolution $d \omega = 0.05$ 
for each different color and the best is chosen. In section~\ref{Locallyadaptivedamping}, a downhill 
simplex method~\cite{Nel-1965} is used that turned out to be stable and 
additionally much more computationally efficient. Both methods were compared and achieved very similar 
convergence rates.

\subsection{Regular tetrahedron}
We first turn our attention to the $\omega-$Jacobi,
the lexicographic Gauss-Seidel and the four-color smoother for a regular tetrahedron, that is
characterized by having edges of equal length. In the case of
Jacobi relaxation, we use the optimal parameter $\omega=0.8$ as predicted by LFA. For the Gauss-Seidel and 
4-color smoothers no damping is applied. 
In Table \ref{table-point-1} we show for different pre-smoothing $(\nu_1)$ and post-smoothing $(\nu_2)$ steps
the smoothing factor $\mu$, the two-grid convergence factor $\rho$ as 
predicted by LFA, and the experimental observed
convergence factor obtained by using $W-$cycles. 
Clearly the two-grid convergence factors are very well predicted by LFA
in all cases. We also observe that the four-color smoother provides the best convergence factors among
the considered point-wise smoothers.

It should be noted that the smoothing factors on tetrahedra are somewhat worse than one would expect for triangles in 2D. For example, the smoothing factor of Gauss-Seidel on an equilateral triangle is 0.416 in 2D, while it is 0.521 for a regular tetrahedron in 3D. This behaviour is similar to what is observed for standard finite difference discretizations of the 3D Poisson equation (see, for example \cite{TOS01}). 

In passing, note that it is possible to improve the convergence factors of the four-color smoother by using relaxation parameters. For instance, in the case of  $(\nu_1,\nu_2)=(2,1)$, the convergence factor for four-color smoother can be improved from $0.153$ to $0.090$ by using overrelaxtion parameters $\omega_i = 1.15, 1.20, 1.25, 1.25$ as introduced in (\ref{operatorS}). Similarly for lexicographic Gauss-Seidel, the convergence factor can be improved from $0.176$ to $0.141$ with the relaxation parameter $\omega=1.2$.

\subsection{Optimized tetrahedron}
\begin{figure}
  \caption{\label{fig:optimizedTet} Optimized tetrahedron with two non-unit edge lengths}

  \vspace*{1.0ex}
  \centering
  \includegraphics[width=0.4\textwidth]{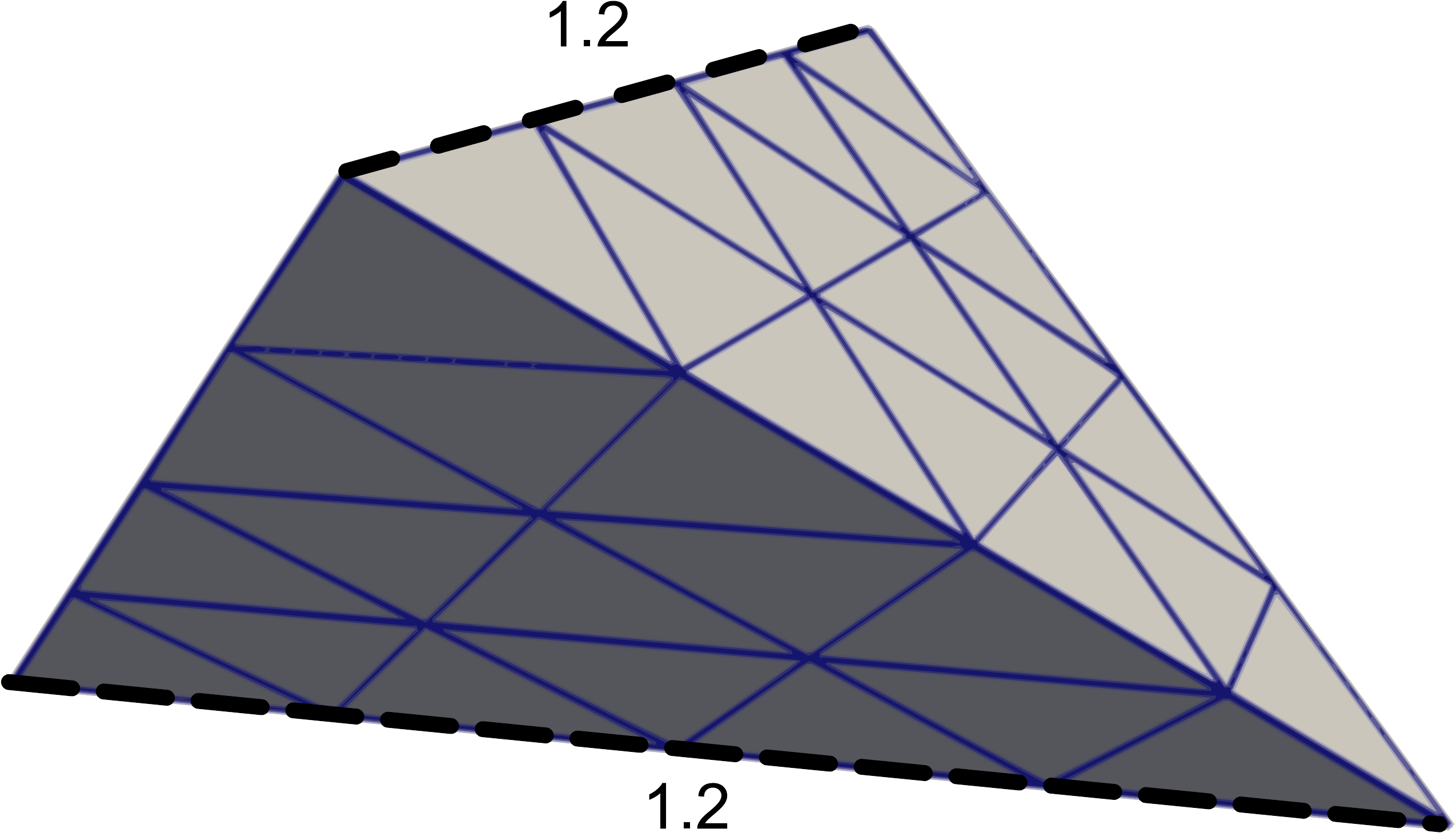}
\end{figure}

A regular input tetrahedron is not the best shape for our regular refinement with respect to the convergence rate. Inside a structured region of a refined coarse grid element, one pair of stencil entries relates to an edge of the element. After one refinement step, four child elements are similar to the parent element, only translated and scaled by the factor $0.5$. The remaining four child elements contain six edges of its parent element, which are translated and scaled by the factor $0.5$. The remaining seventh edge is introduced by Bey's refinement \cite{Bey95}. This edge is located between the midpoint of edge $e_{02}$ and the midpoint of edge $e_{13}$ (see Figure~\ref{regtet}) of the parent element. The subscripts denote the vertex indices of the refined element. In the case of a regular tetrahedron, all edges except of the additional introduced edge have the same length. In contrast to the six other entry pairs, which are typically negative, the seventh pair has positive entries in the stencil. Thus, compared to the other edges there is a weak connection, if this edge is small. For example by setting the edge lengths of $e_{02}$ and $e_{13}$ equal to $1.2$ and keeping all other four edge sizes equal to $1.0$, we receive an optimized tetrahedron concerning the convergence rate (see Figure~\ref{fig:optimizedTet}).

The exceptional convergence factors associated with the four-color smoother for this optimized tetrahedron are pointed out in Table \ref{table-point-2}. On this type of tetrahedra, the convergence factor
of the four-color smoother is approximately $\rho=0.106$ for one pre- and post-smoothing steps,
having at the same time the desirable property of high parallelism. This factor can also
be further enhanced to $0.080$ using damping parameters of values $\omega = (1, 1, 1, 1.2)$.

\begin{table}[h]
\begin{center}
\begin{tabular}{|c|c|c|c|c|c|c|c|c|c|}\cline{2-10}
\multicolumn{1}{c}{} & \multicolumn{3}{|c|}{Damped Jacobi} &
\multicolumn{3}{|c|}{Gauss-Seidel} & \multicolumn{3}{|c|}{Four-color}\\
\hline $\nu_1, \nu_2$ & $\mu^{\nu_1+\nu_2}$ & $\rho$ & $\rho_h$ & $\mu^{\nu_1+\nu_2}$ & $\rho$ & $\rho_h$ & $\mu^{\nu_1+\nu_2}$ & $\rho$ & $\rho_h$\\
\hline $1,0$ & 0.741 & 0.640 & 0.637 & 0.521 & 0.434 & 0.427 & 0.500 & 0.407 & 0.389 \\
$1,1$ & 0.550 & 0.409 & 0.407 & 0.272 & 0.223 & 0.219  & 0.250 & 0.195 & 0.197 \\
$2,1$ & 0.406 & 0.298 & 0.297 & 0.141 & 0.176 & 0.174 & 0.125 & 0.153 &  0.153 \\
$2,2$ & 0.301 & 0.250 & 0.250 & 0.074 & 0.143 & 0.141 & 0.062 & 0.123 & 0.125\\ \hline
\end{tabular}
\end{center}
\caption{LFA smoothing factors, $\mu^{\nu_1+\nu_2}$, LFA
two--grid  convergence factors, $\rho$, and  measured $W$--cycle convergence
rates, $\rho_h$, for a regular tetrahedron} \label{table-point-1}
\end{table}

\begin{table}[h]
\begin{center}
\begin{tabular}{|c|c|c|c|c|c|c|c|c|c|}\cline{2-10}
\multicolumn{1}{c}{} & \multicolumn{3}{|c|}{Damped Jacobi} &
\multicolumn{3}{|c|}{Gauss-Seidel} & \multicolumn{3}{|c|}{Four-color}\\
\hline $\nu_1, \nu_2$ & $\mu^{\nu_1+\nu_2}$ & $\rho$ & $\rho_h$ & $\mu^{\nu_1+\nu_2}$ & $\rho$ & $\rho_h$ & $\mu^{\nu_1+\nu_2}$ & $\rho$ & $\rho_h$\\
\hline
$1,0$ & 0.720 & 0.602 & 0.598 & 0.492 & 0.401 & 0.392 & 0.442 & 0.345 & 0.331 \\
$1,1$ & 0.517 & 0.362 & 0.360 & 0.243 & 0.151 & 0.145 & 0.196 & 0.106 & 0.105 \\
\hline
\end{tabular}
\end{center}
\caption{LFA smoothing factors, $\mu^{\nu_1+\nu_2}$, LFA
two--grid  convergence factors, $\rho$, and  measured $W$--cycle convergence
rates, $\rho_h$, for an optimized tetrahedron} \label{table-point-2}
\end{table}

\begin{figure}[!htb]
\subfigure[Needle]{\includegraphics[width=0.25\textwidth]{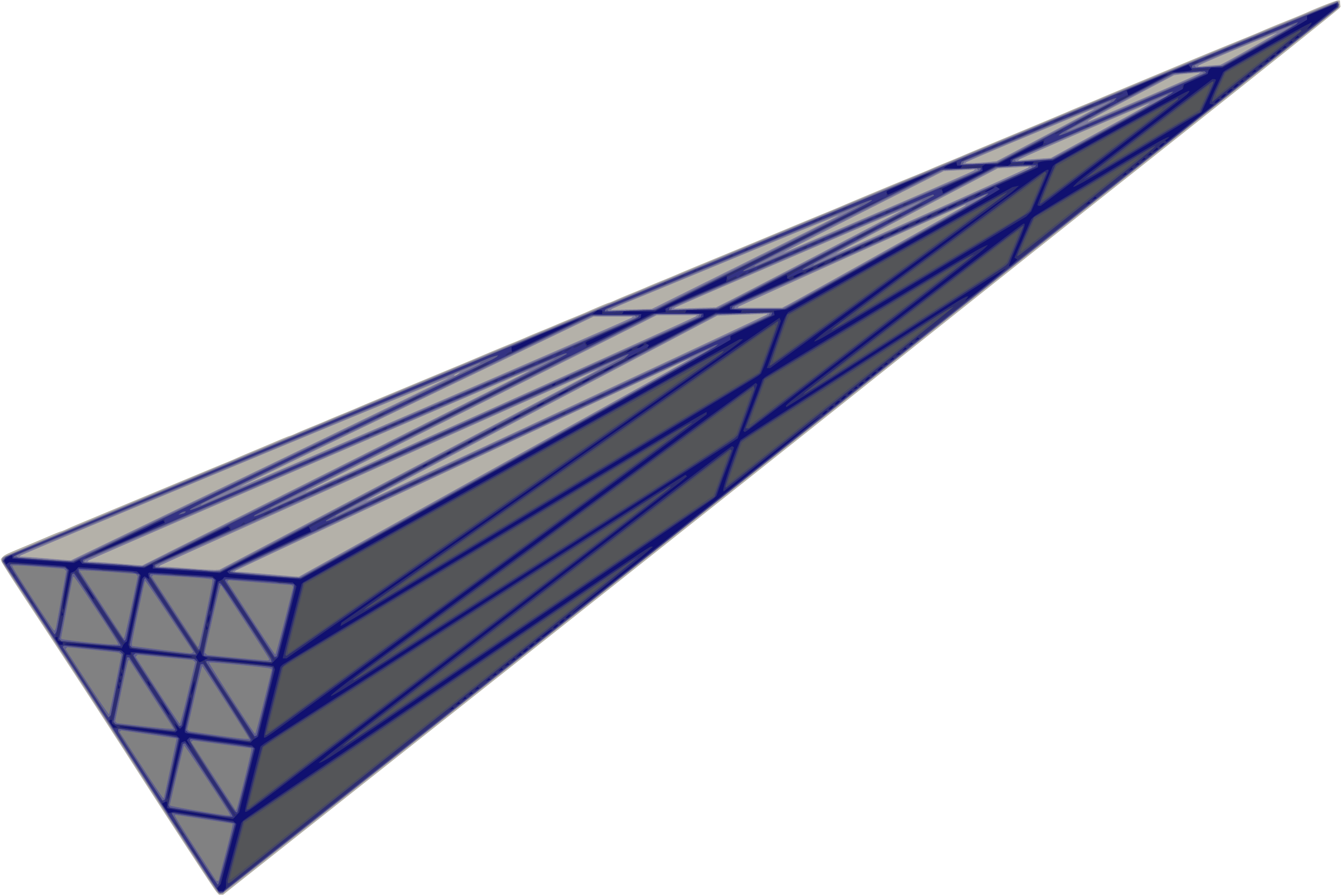}}\hfill
\subfigure[Wedge]{\includegraphics[width=0.25\textwidth]{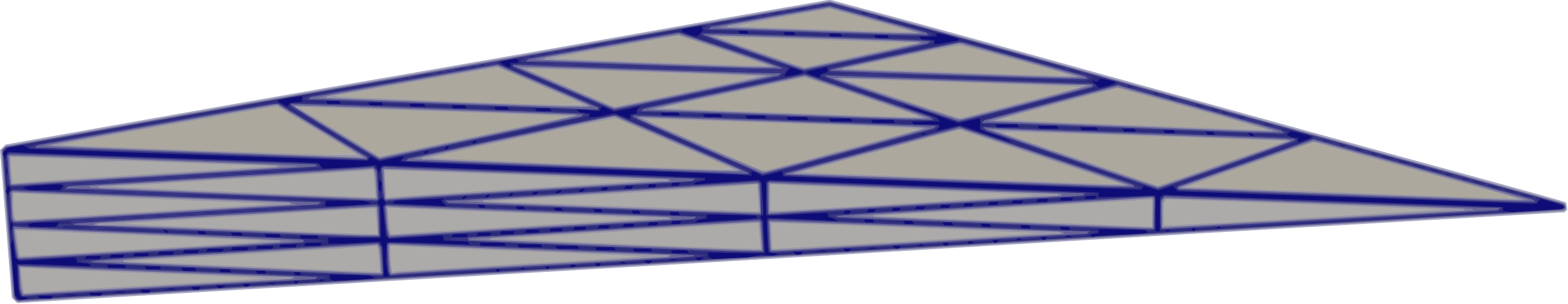}}\hfill
\subfigure[Spindle]{\includegraphics[width=0.20\textwidth]{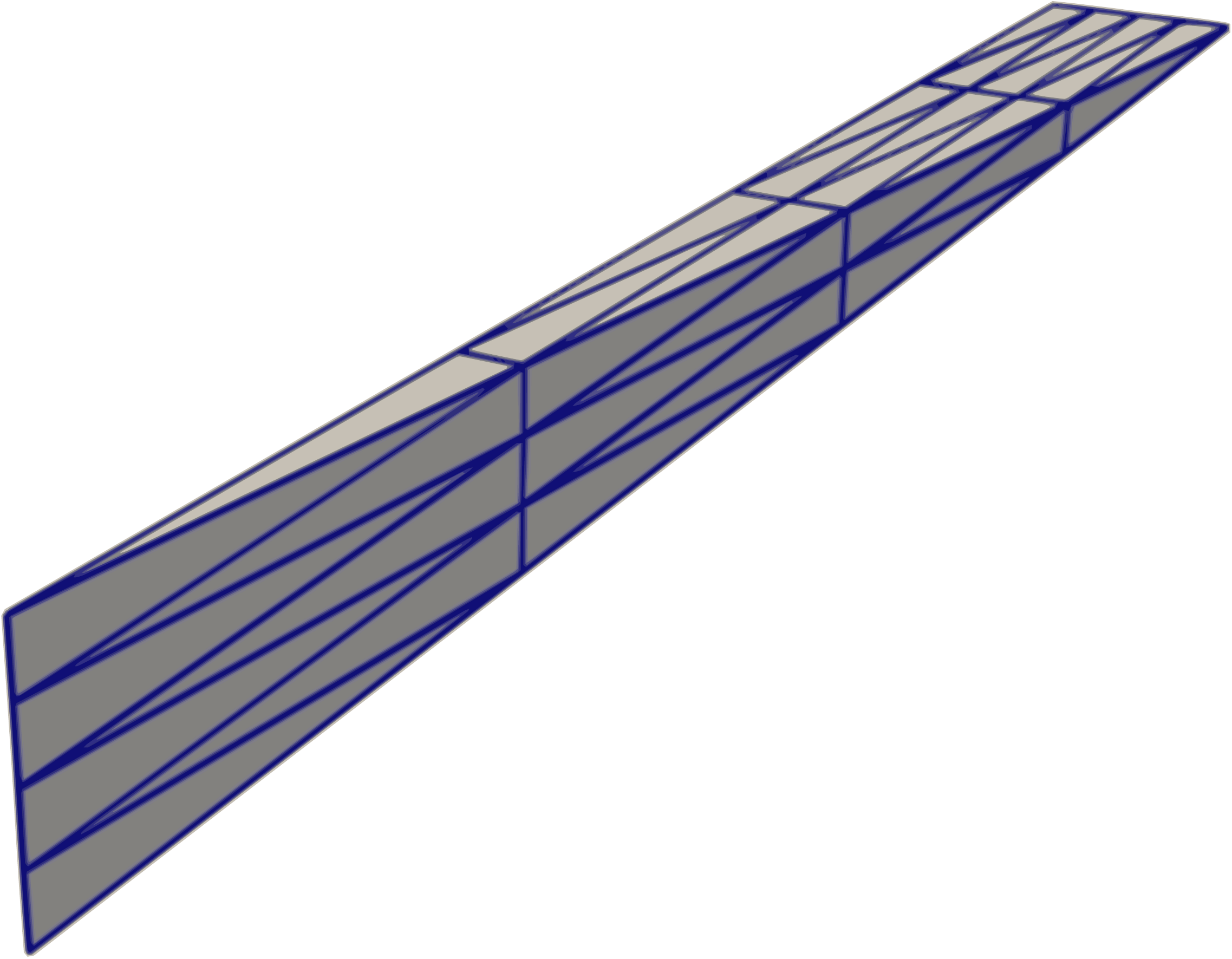}}\\
\subfigure[Spade]{\includegraphics[width=0.25\textwidth]{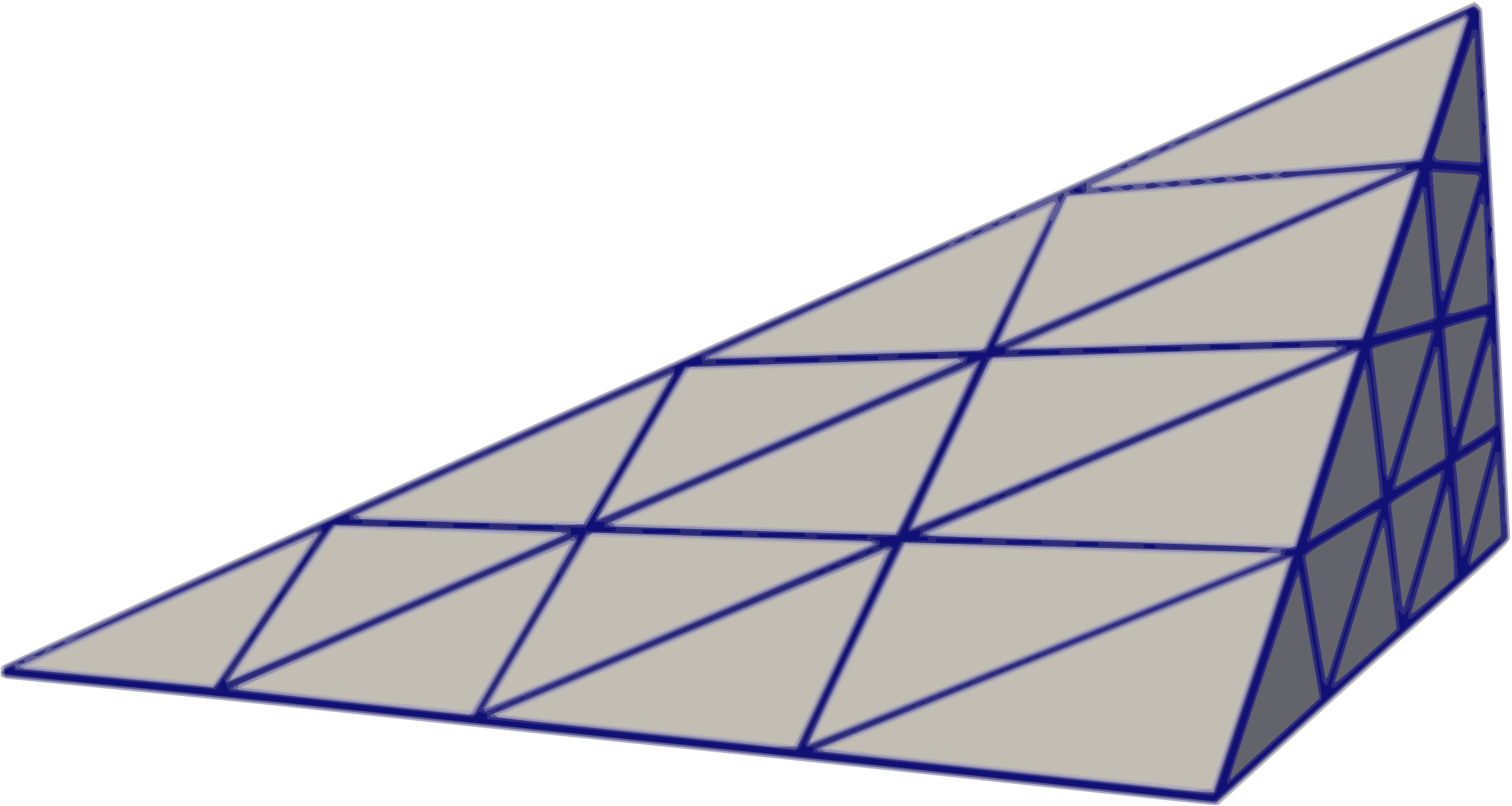}}\hfill
\subfigure[Sliver]{\includegraphics[width=0.3\textwidth]{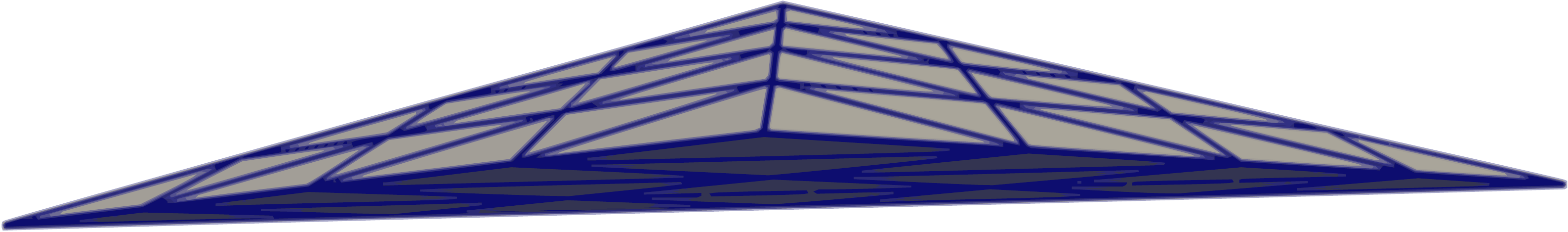}}\hfill
\subfigure[Cap] {\includegraphics[width=0.25\textwidth]{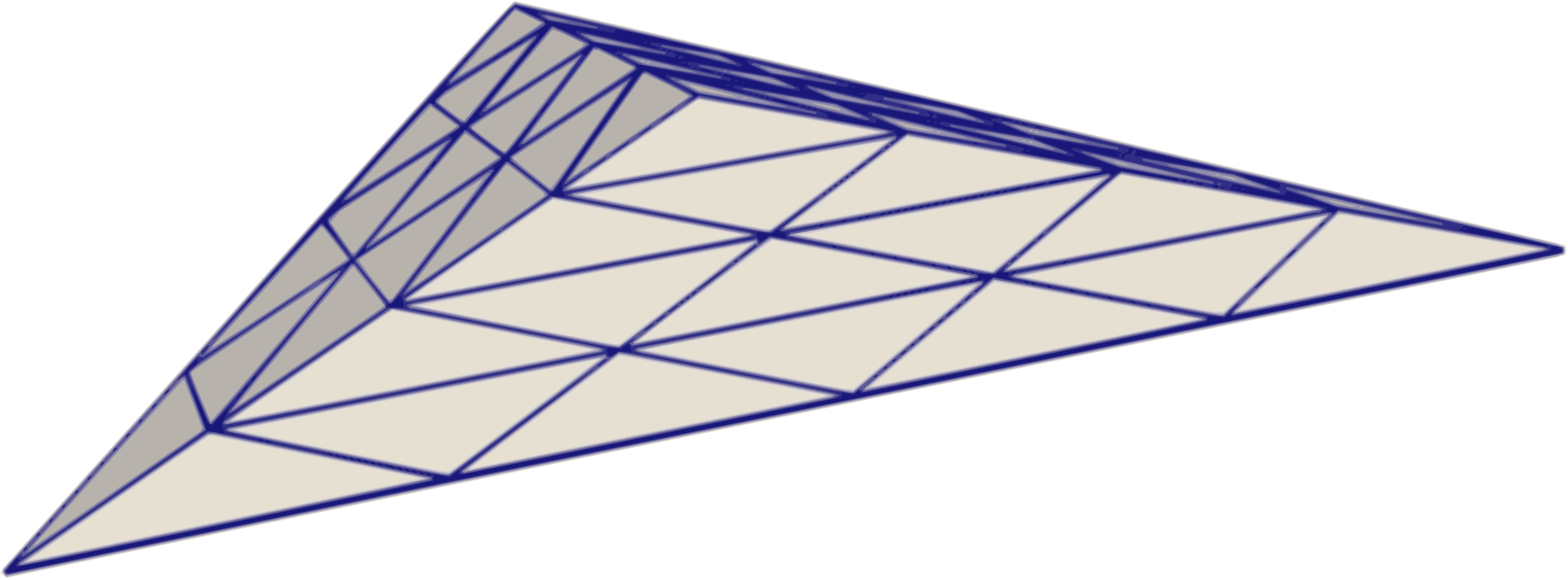}}\\
\caption{Different classes of degenerated elements}
\label{fig:bad-tetrhedron}
\end{figure}

\subsection{Degenerated tetrahedra}
On the other hand, the highly satisfactory convergence factors obtained for regular and
optimized tetrahedra, can deteriorate quickly when poorly shaped tetrahedra are treated. In two dimensions,
there are only two types of failure, angles close to $0^0$ and angles close to $180^0$, and no failures
of the first kind implies no failures of the second. In three dimensions, we can classify poorly
shaped tetrahedra according to both dihedral and solid angles. Dihedral is the interior angle between two face planes, 
while solid angles are formed by the intersection of three planes at one vertex. There are then six types
of bad tetrahedra, as shown in Figure~\ref{fig:bad-tetrhedron}. 
%A needle or thin has one small solid angle, but not small or large
%dihedrals. A wedge-like element has small but not large dihedrals and no large angles of any kind. A spindle
%allows large but not small dihedrals, and small but not large solid angles. A spade element has one large solid angle as well as small dihedral angles. A sliver has small and large dihedrals, but no large solid angles.
%A cap-like tetrahedron has a large solid angle. 

Table \ref{bad} shows the characteristics of differently shaped tetrahedra. Besides the number of solid and dihedral angles, the edge ratio, minimal and maximal angles between two edges are provided for an exemplary tetrahedon of this type. 
For these elements, the smoothing factor and the two-grid convergence factor are given for
one pre- and one post-smoothing step ($\nu_1=1, \nu_2=1$). Poor behavior of the
four-color smoother for these type of tetrahedra is shown together with the good performance that is observed for the
regular and optimized tetrahedra.

\begin{table}[h]
\begin{center}
\scriptsize 
\begin{tabular}{|c|c|c|c|c|c|c|c|c|}\cline{2-9}
\multicolumn{1}{c|}{}
                      & Regular & Optimized & Needle   &  Wedge  & Spindle & Spade   & Sliver  & Cap \\
\hline
  Small solids & - & - & 1 & - & - & - & - & - \\
  Large solids & - & - & - & 2 & - & 1 & - & 1 \\
\hline
  Small dihedrals & - & - & - & 1 & 2 & 2 & 4 & 3 \\
  Large dihedrals & - & - & - & - & - & 1 & 2 & 3 \\
  \hline
  $\mu^{\nu_1+\nu_2}$ & $0.250$ & $0.196$   & $0.982$  & $0.786$ & $0.980$ & $0.590$ & $0.896$ & $0.872$    \\
  $\rho$              & $0.195$ & $0.106$   & $0.982$  & $0.740$ & $0.980$ & $0.500$ & $0.889$ & $0.739$    \\
  \hline
  Edge ratio & $1.0$ & $1.15$   & $10.0$   & $4.0$    & $10.0$   & $1.67$    & $1.40$ & $1.71$  \\
  Min angle & $60^0$ & $54.9^0$ & $5.7^0$  & $14.3^0$ & $5.7^0$  & $33.6^0$  & $45.6^0$ & $31.4^0$  \\
  Max angle & $60^0$ & $70.2^0$ & $87.1^0$ & $82.8^0$ & $87.1^0$ & $112.9^0$ & $88.9^0$ & $117.2^0$ \\
  \hline
  
\end{tabular}
\end{center}
\caption{LFA smoothing factors, $\mu^{\nu_1+\nu_2}$, and two--grid  convergence factors, $\rho$, of the four-color smoother
for different tetrahedra} \label{bad}
\end{table}

\subsubsection{Wedge}
To overcome this difficulty, we propose specific smoothers for some of the poorly
shaped tetrahedra. The first type of such tetrahedra is the wedge type, which has
one edge much smaller than the others. This causes two entries of the stencil
to be much larger than the others, producing the strongest connection in one direction.
Since the common lore states that the good smoothing of errors is obtained when the strongly
coupled unknowns are collectively relaxed, a line smoothing in the direction of the
smaller edge is expected to improve the performance. In fact, a two-grid convergence factor $0.122$ is
obtained by using the corresponding line smoother with one pre- and one post-smoothing steps.
For very anisotropic wedge-type tetrahedra improved convergence rate easily compensated the 
somewhat higher computational cost of the tridiagonal solvers in the line smoother.
Table \ref{convdiffsmooth} presents an overview of the effect of line-wise and plane-wise 
smoothers onto the convergence rate for this and the following tetrahedral shapes.
\begin{table}[h]
\begin{center}

\begin{tabular}{|c|c|c|c|c|c|}\cline{2-6}
\multicolumn{1}{c|}{}
                          &  Wedge  & Spade & Needle   & Spindle  & Sliver/Cap \\
\hline
  Line-wise & good & middle & poor & poor & poor \\
  Plane-wise & good & good & good & good & poor  \\
\hline
\end{tabular}
\end{center}
\caption{Smoothing properties of line-wise and plane-wise smoothers
for different tetrahedra} \label{convdiffsmooth}
\end{table}

\subsubsection{Spade}
Another tetrahedral shape is generated by assuming two short edges, which are connected to each other.
Keeping the other edge lengths similar but larger results in a spade. In comparison to the needle-type, discussed next, only two stencil-pairs of three in a plane are large. Already an improvement is visible, choosing a line-wise smoothing in the direction of one short edge. With two smoothing steps in total, the computational work is comparable to point-wise smoothing, but results in a better convergence rate of $0.23$. Hereby it yields to similar convergence, regardless of smoothing alternatively in both directions of the short edges, or repeatedly along one of them. Point-wise smoothing only provides a convergence rate of $0.59$ (see Table~\ref{bad}). 

Another possibility is applying a zebra-plane smoothing about the faces spanned by those edges. With only one smoothing step, a plane-zebra smoother yields convergence factors about $0.105$. Whether or not this reduces the overall cost depends on the implementation of the plane-smoother and the degree of anisotropy. Indeed, because of significant higher costs of plane-wise smoothing in our current implementation, discussed in the next section, we would choose line-wise smoothing for that example.

\subsubsection{Needle}
The next bad tetrahedron is the needle-type, which is
characterized by the fact that one vertex is far away from the
opposing face. This feature produces a stencil with stronger connections in
the direction of the plane defined by the face opposite to such vertex. Therefore,
the corresponding plane relaxation will result in excellent smoothing. In fact, a very
good convergence factor about $0.121$ is reached with only one smoothing step for a zebra-plane
smoother. We want to emphasize the fact of the use of zebra-plane smoothers rather than lexicographic
plane relaxation. This is due to the significant improvement of the convergence factors when zebra-type
smoothers are used. For example, in this case a factor of $0.330$ is provided by lexicographic plane relaxation.
Thus one plane-zebra sweep is roughly as effective as two lexicographic plane sweeps $((0.330)^2 = 0.1089)$. 
One might think that alternating line relaxation with respect to the two directions defining
the corresponding face could be a suitable choice, but very poor convergence factors are provided
by this strategy. For instance, only a factor of about $0.942$ is obtained when using
two smoothing steps. We remark that all strongly connected unknowns
should be relaxed collectively and therefore a plane smoother must be used.

\subsubsection{Spindle}
Analogously, the behavior of the spindle-type is similar to that of needle type. In this case
the strong connection also corresponds to the direction of a face of the tetrahedron, and therefore
a plane smoothing will be a good choice. Again, very satisfactory convergence properties are provided
by using such smoother. LFA predicts a two-grid convergence factor about $0.124$ with only one
smoothing step if the adequate zebra-plane smoother is activated, that is, the zebra-plane relaxation 
updating all strongly connected unknowns.

\subsubsection{Sliver and cap}
According to Shewchuk \cite{Shew02} "good" tetrahedra are of two 
types: those that are not flat, and those that can grow arbitrarily flat without having a large planar 
or dihedral angle. The "bad" tetrahedra result in error bounds that explode, and a dihedral angle or a planar angle
that approaches $180^0$, as they are flattened. Sliver as well as cap type can grow arbitrary flat and have in addition 
large planar or dihedral angles and thus have bad interpolation properties.

For a different reason, we also were not able find an optimal collective relaxation method for sliver and cap type. 
Here, the strongly connected points are not located on an edge or plane; therefore the smallest block would already include all grid points.
In principle alternating plane relaxations, or alternating line relaxations 
improve the convergence rate compared to a point-wise relaxation like shown 
for a spade. However, similarly in this case the convergence rate deteriorates when the degree 
of anisotropy increases.

However, this is not a limitation from the practical point of view. Many Delaunay mesh generators
include a post-processing to find and remove sliver and cap type tetrahedra, because of the poor approximation quality
of the finite element methods using these types.

%% file: results.tex
In this section we will discuss numerical experiments to illustrate  two different strategies, 
from subsection~\ref{sec:smoothersemistructuredmeshes} and \ref{sec:localadapt}. Firstly, we
will demonstrate how to use LFA to design 
efficient multigrid solvers for semi-structured tetrahedral grids. We will use the 
block-wise processing to choose different smoothers depending on the
shape of the tetrahedra. This kind of strategy has already been used successfully for
two-dimensional linear-elasticity problems in \cite{Gaspar:2010:NLAA}. Here, we will apply
it for three-dimensional tetrahedral meshes. The second strategy consists of applying LFA to improve the 
overall convergence rates by adapting smoothing parameters. This include the damping parameters, 
number of pre- and post-smoothing steps, and smoother types individually for each block, in a complex domain.

\subsection{\label{sec:smoothingCosts}Computational cost for our implementation}

Especially for the more sophisticated smoothers it is important to keep their computational costs in mind. The computational costs of an algorithm is often expressed in arithmetic operations. An advantage of this measure is its platform independence. Unfortunately, it often cannot give realistic predictions of performance for an implementation on specific hardware. The runtime on many modern architectures is dictated by memory bandwidth, cache effects, instruction level parallel execution, or memory alignment issues. 
Here we will present an overview of the smoothers that are available in HHG. To allow a fair comparison, all algorithms are implemented in an efficient way, but are not yet optimized specifically for each architecture. The performance is evaluated on one core of a \emph{Xeon 5550 Nehalem} chip with 8 MB shared cache per chip. The code was produced by the \emph{Intel 10.1} compiler (\emph{flags: -O3}). We will measure the required number of clock cycles per unknown for three different block sizes for a structured region: $65$, $129$ and $257$ unknowns in each direction (see Table \ref{compCosts}). Most differences can be attributed to cache-effects. Applying damping parameters cost additional $3$ clock cycles or less. This is not really severe for any of the schemes introduced. Thus it should be applied, if it leads to advantages. 

The compiler is basically able to utilize SSE SIMD instructions (small vector processing) for the stencil evaluations. For large vectors, a 4-color smoother with a splitting of the colors to different arrays would be necessary to achieve a good performance. Further, the data layout on the structured regions would also fit well to a parallel execution on GPUs. Here, especially smoothers with a low degree of dependency (like a 4-color smoother, zebra-line smoother, or a parallel plane-wise smoother) would be most suitable.

Especially to perform line-wise or plane-wise smoothing over interface boundaries, using a multigrid in one or two dimensions seems to be a promising approach. However it not always obvious to define a line or plane for the whole domain, since the coarsest mesh is unstructured.

\paragraph{Line-wise smoother}
The line-wise smoother is implemented by the Thomas' Algorithm or tridiagonal matrix algorithm, a simplified version of Gaussian elimination. For constant coefficients, terms involving only the coefficients of the stencil are calculated once per stencil. Thus the required operations per unknown are comparable to a Gauss-Seidel smoothing. The number of clock cycles is only slightly higher. 
%in x-direction and up to two times higher for other directions. The reason is that in our implementation the continuous memory is aligned to x-direction. 
A zebra-line fashion has a high degree of parallelism in two directions. A parallel execution in the line-wise direction can be achieved by cyclic reduction (see e.g. \cite{Hoc-1965}). 

\paragraph{Plane-wise smoother}
The advantage of an iterative method for a plane-wise smoother is that the plane does not have to be solved exactly. A good approximation is enough to obtain a optimal convergence rate of the whole V-cycle. This behavior was observed e.g. in \cite{Llo-2000}, where a 2-dimensional multigrid is used for each plane.  

For tetrahedral structured regions, it is not straightforward to apply plane-wise multigrid, since each plane sizes has a different size. 
Therefore, we solved arising system by a Conjugate Gradient (CG) algorithm. Since the number of CG iterations increases with respect to a growing plane-size, also the number arithmetic operations and thus solving time are effected.

\begin{table}[h]
\begin{center}
\begin{tabular}{|l|c|r|r|r|}\cline{1-5}
\multicolumn{1}{|c|}{} & \multicolumn{1}{|r|}{} &
\multicolumn{3}{|c|}{\,\, Clock cycles (line sizes) } \\
%\hline
\cline{3-5}
Algorithm & Operations & 65 & 129 & 257 \\
\hline
lex. Gauss-Seidel       &  29 & \,\,\,\,\,\,\,\, 35 & \,\,\,\,\,\,\,\, 35 & \,\,\,\,\,\,\,\, 36 \\
Four-color         &  29 &  40 & 39 & 40 \\
Line-wise        &  29 &  40 & 39 & 41 \\
%non x-line-wise    &  29 &  48 & 53 & 81 \\
Plane-wise      & 270-600 &  360 & 456 & 694 \\
%cg                 & 305- &  485 & 809 & 1626 \\
%xy-Line-wise   & 29 & 48-88 & 33 & 48-88 \\
\hline
\end{tabular}
\end{center}
\caption{Arithmetic costs per unknown and measured clock cycles for a Gauss-Seidel, four-color, line-wise and plane-wise smoothers} \label{compCosts}
\end{table}

\subsection{\label{sec:smoothersemistructuredmeshes}Smoothers on semi-structured meshes}

In the previous section we presented the impact of strongly
anisotropic meshes on the convergence rate for structured tetrahedral meshes.
Furthermore we found appropriate line- and plane-smoothing strategies for different
types of element degeneration. In a smoothing step of HHG, we first smooth the  
the interface points (macro vertices, edges, and faces) between all structured regions. 
Afterward we smooth the structured regions themselves.
While for point-wise smoothing this proceeding has nearly no influence on the
convergence \cite{Huelsemann:2005}, we will study its influence for plane- or
line-smoothing in the first numerical experiment, where we measure
the impact of the interface points between the structured regions using semi-structured
meshes. In the second subsection, the advantages of using different smoothers for different
shaped tetrahedra composing the domain are displayed.

\subsubsection{\label{sec:sameanisotetr}Anisotropic smoothing}

\begin{figure}[!htb]
\caption{Connecting faces in different directions}
\subfigure[Connection in isotropic direction]{\includegraphics[width=0.5\textwidth]{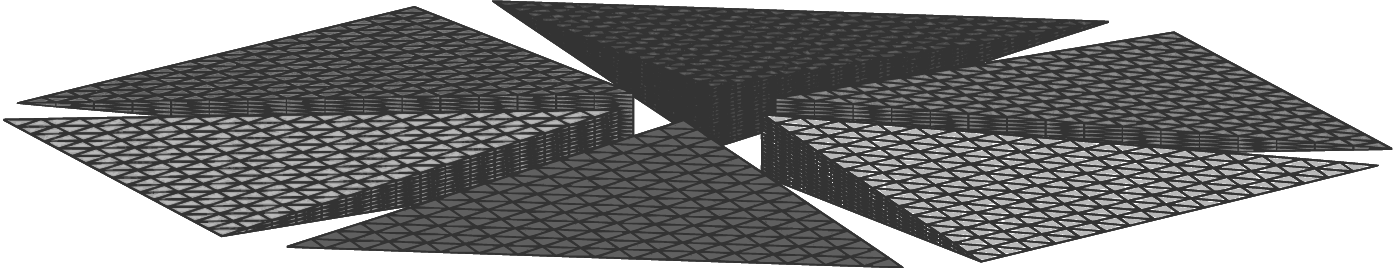}}\hfill
\subfigure[Connection in anisotropic direction]{\includegraphics[width=0.45\textwidth]{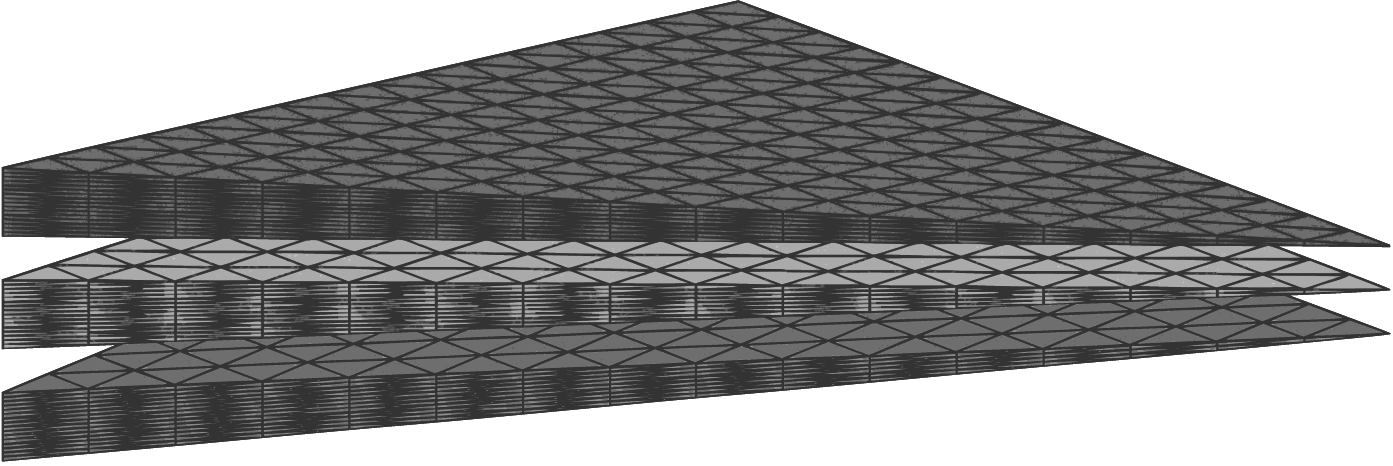}}\hfill
\label{fig:wedge_conn}
\end{figure}

In the following two numerical experiments we construct our computational domain
out of the wedge type elements. Figure~\ref{fig:wedge_conn}~(a) shows a setup
where six structured regions are sharing anisotropic faces. In Figure~
\ref{fig:wedge_conn}~(b) the same type of tetrahedra are used, but are here
connected at their regular triangular face. In the following part, we apply two
smoothing steps for all smoothers. The convergence factors are measured after 50
W--cycles. We point out that the used elements have a edge ratio of $10$
between the shortest and the largest edge length and dihedral angles of $5.7^0$.
Otherwise the following discussed effects are not so severe. On the fine grid
we end up with $1.73 \cdot 10^7$ elements in case (a) and $3.46 \cdot 10^7$ elements
in case (b).

For a single region, lexicographic Gauss-Seidel achieves a convergence
factor of $0.934$. This factor is nearly the same for both domains (a)
and (b). In contrast a line smoother yields a convergence factor of
$0.072$ for a single region. For domain (a) the interfaces slightly
effect the convergence rate to $0.094$. However it is worse if we cannot
solve the full line of unknowns along the strong coupling, but only parts of the
line. This is case for domain (b) since two interfaces are breaking these lines
and do not allow a fully collective line-smoothing. This explains a bad
convergence factor of $0.847$ for (b). Damping parameter of $1.05$ inside the
structured regions and $1.15$ for the boundary points help a bit to
decrease the convergence factor to $0.796$. This problem has been reported
in some papers about grid-partitioning, see for example \cite{oos}. 
If possible, this case should be avoided in the original partitioning of the domain.

The same behavior is observed by applying plane smoothing to adjoined needle or
spindle type regions. In these cases lines of unknowns, which are not
collectively updated since they are interface points, increase the convergence
rate.

\subsubsection{\label{sec:localsmoothers}Locally adaptive smoothing}

\begin{figure}
  \caption{\label{fig:diffTetrDomain} Domain consisting of wedge, regular and needle tetrahedral types}
  \vspace*{1.0ex}
  \centering
  \includegraphics[width=\textwidth]{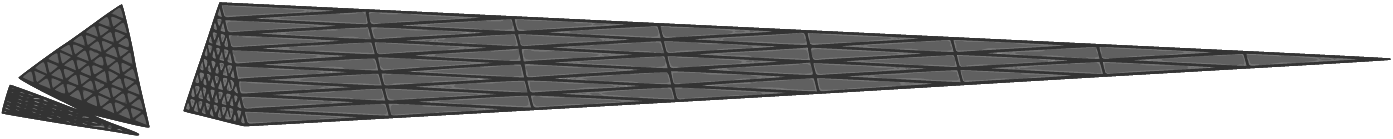}
\end{figure}

Contrary to the previous examples, we apply different smoothers according to the
geometry of the structured region. A model domain consisting of four different
tetrahedral types is presented in Figure~\ref{fig:diffTetrDomain}. Wedge and
needle elements share some of their faces with regular elements. The whole domain
consists of $1.73 \cdot 10^7$ tetrahedron elements on the finest level. Table
\ref{diffshapes} shows some properties of the single structured regions without its
interfaces to any other regions. \\

In order to get a similar convergence rate for the composed domain, we choose different smoothers
according to the underlying mesh geometry. The regular tetrahedra are smoothed by using four-color 
relaxation with damping parameters $\omega = (1.15, 1.20, 1.25, 1.25)$, and 
two pre- and one post-smoothing steps. From \Sec{sec:LFA_results}, we know that a convergence factor of about
$0.09$ is achieved for these regular tetrahedra. For the wedge element, again a value of $0.09$ is predicted
by LFA, if a line-smoother is used with two pre- and one-smoothing steps. By the other hand, for the needle 
tetrahedron a zebra-plane smoother with only one smoothing step is enough to reach a convergence factor
about $0.12$. Taking into account these local convergence factors and considering the worst of them, 
we can predict a global factor of $0.12$ for the composed domain. \\
 
By performing the multigrid algorithm with the corresponding local components chosen above, 
a convergence factor of $0.17$ is obtained for the whole computational domain. An additional
smoothing step of the macro faces after smoothing the interior was applied. Using damping
parameter at the interface points and in the structured regions, the convergence
factor improves to $0.14$, which is very close to the expected value of $0.12$. As a
reference, a standard lexicographic Gauss-Seidel with $\nu_1, \nu_2 = 2,1$
provides a convergence factor of $0.95$. It turns out that using different smoothing
strategies depending on the sub-domains is very advantageous and following this rule, very efficient multigrid
solvers can be designed for quite complicated three-dimensional domains. However to be fair,
one has to consider the additional computational overhead for line- or
plane-smoothing. If and how much it pays out can also depend much on the
implementation of the smoothers.

\begin{table}[h]
\begin{center}
\begin{tabular}{|l|c|c|c|c|c|c|}
\hline
\multicolumn{1}{|r|}{Shape} & edge ratio & max angle & min angle & smoothing & $\nu_1, \nu_2$ & $\rho$ \\
\hline
Needle & 10 & $87^0$ & $5.7^0$ & plane-wise & 1,0 & 0.12 \\
Regular & 1 & $60^0$ & $60^0$ & 4-color & 2,1 & 0.09 \\
Wedge & 5 & $110^0$ & $11^0$ & line-wise & 2,1 & 0.09 \\
\hline
\end{tabular}
\end{center}
\caption{Convergence factor for the single structured regions} \label{diffshapes}
\end{table}

\subsection{\label{sec:localadapt}Locally adaptive damping and smoothing steps on semi-structured meshes}

In order to achieve a good overall convergence of the solver, we try to
distribute the computational work such that all regions of the domain receive a
similar error reduction in each smoothing step. In our approach we use LFA to control 
this process. Partial smoothing has already been promoted, \eg{}, in \cite{Brandt:1977:MoC}.
In addition to that, we apply individual damping parameters to each structured region.

With HHG's semi-structured meshes it is relatively easy to locally adapt
smoothing parameters to mesh-induced convergence rate variations. Since each
coarse mesh element gets refined in a regular way (\cf{} \Fig{fig:2d_levels}),
all its sub-elements are similar, and the operator stencils of all interior grid
points are the same. Therefore, one can determine the optimal smoothing
parameters for each coarse grid element and use them for all its
sub-elements. We will use LFA prediction to optimize the convergence rate of the
structured regions by choosing suitable the smoother parameters. Further we will 
provide suggestions how to deal with interface points between those regions.  
Step by step we 
will apply different modifications to the smoothing procedure and observe its 
influence. Thus we end up with an asynchronous or chaotic smoothing \cite{Barlow:1982,Baudet:1978,Chazan:1969}. 
However, this is usually a result of an efficient parallel implementation anyway.

% The look-up table approach first classifies tetrahedral elements by their
% shape. We use the classification into nine tetrahedron types according to edge
% length ratios. Two examples are the needle and the wedge shown in
% \Sec{sec:LFA_results}. In an off-line step the convergence rates of each type
% for different degrees of anisotropy are measured and stored in a look-up
% table. As a measure of the anisotropy the aspect ratio of the longest and the
% shortest edge is taken. At the start of a productive run HHG maps the
% tetrahedrons of an arbitrary input mesh to an entry of the look-up
% table. Depending on this information, the global number of smoothing steps is
% redistributed among the structured regions, in order to get similar convergence
% rates for each region.

\begin{figure}[!htb]
\caption{\label{fig:halfsphere}Test domain: 293 elements on the coarsest mesh.}
\centering
{\includegraphics[width=0.6\textwidth]{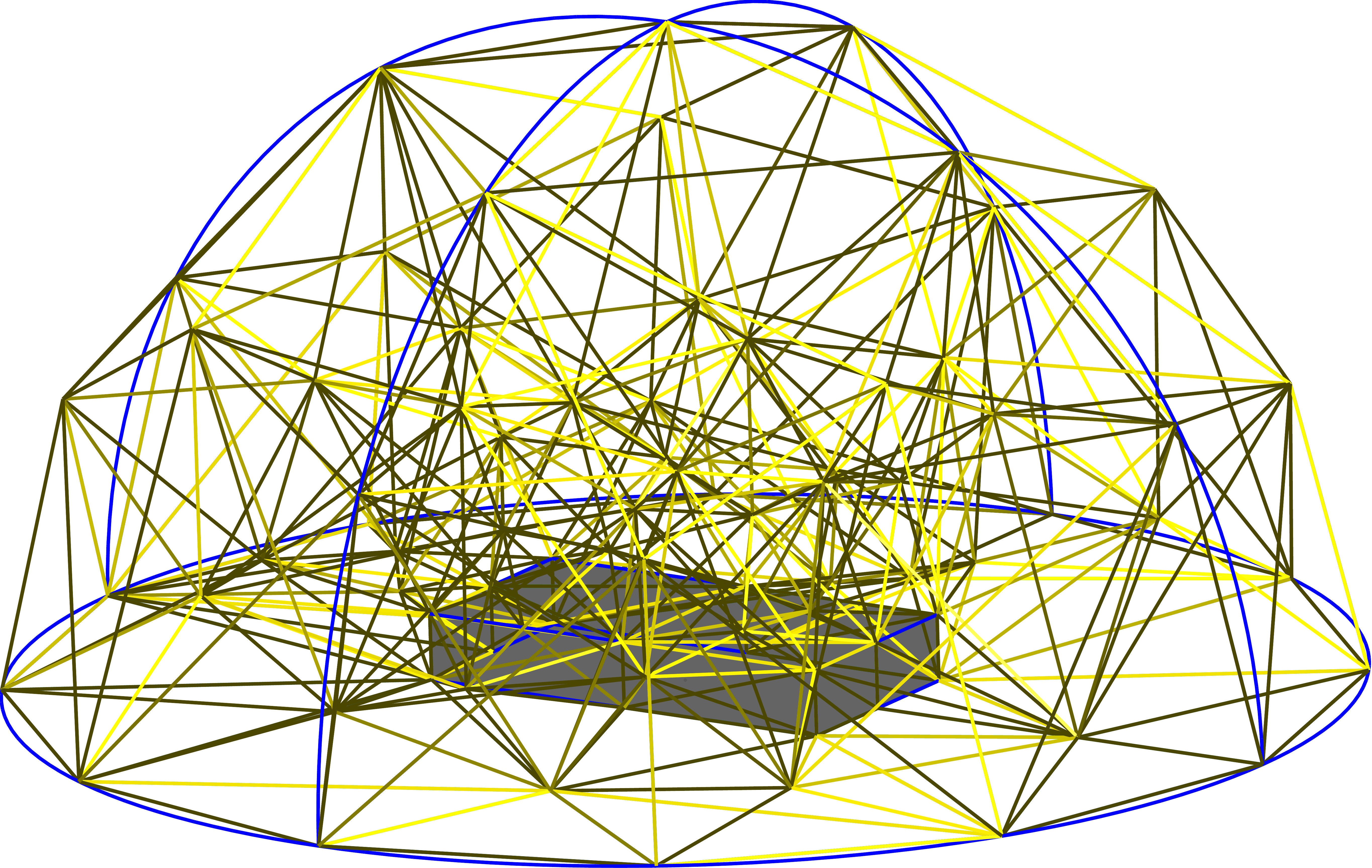}}\hfill
\end{figure}

As a domain for testing the optimizations we use a half ball with a small box
cut out at its bottom. The diameter of the half ball is $6.0$, the box has a
size of $1.0 \times 2.0 \times 0.3$ (see Figure~\ref{fig:halfsphere}). The domain is
discretized with $293$ elements on the coarsest mesh, generated by Gmsh~\cite{Geuzaine:2009:NME}. 

These elements are refined $8$ times, which leads to $1.02 \cdot 10^{8}$ unknowns on the finest
level. A well-known degeneracy measure for FE meshes is the ratio between
inscribed and circumscribed radius $\alpha$, another one is the ratio between
the shortest and the longest edge $\beta$. The average degeneracy measures of
our test mesh are $\alpha = 0.610$ and $\beta = 0.537$. The minimal and maximal
values are $\alpha_{min} = 0.313$, $\alpha_{max} = 0.900$, $\beta_{min} =
0.214$, and $\beta_{max} = 0.830$.

\subsubsection{Unoptimized version and exact block-wise solving}
As a reference a four-color smoother is considered in the domain. For the
unoptimized version $\nu = (4,4)$ yielded the best time to solution. With this
number of smoothing steps, 21 W- or V-cycles are required to reduce the error of a
random initial solution below $10^{-7}$.  
On a single core, one V-cycle takes approx. 16 seconds. For all following experiments, we measure
the convergence rate of the last cycle (see Table~\ref{convHalfSphereResults}). \\

\noindent Further we impose two conditions:

\begin{enumerate}
\item The required amount of communication is the same for all the methods. We are exchanging eight times all ghost boundary points at the interfaces during the smoothing procedure.
\item The total smoothing workload has to be comparable to the unoptimized version in the solving phase. The additional time for the setup phase strongly depends on the implementation of a LFA evaluation in some of the following smoothing strategies. \\
\end{enumerate} 

\noindent The first condition provides an estimate of a tight lower limit for the global convergence rate for our experiments. 
For standard iterative smoothers, we can estimate this limit by 
solving each block exactly (see Table~\ref{convHalfSphereResults} \textit{exact block-wise solving}). 
However this is just to provide a comparison, since it requires 
much more smoothing work than we need for the unoptimized version. The interface points (macro
vertices, edges and faces) between the structured regions are also considered as small blocks 
in themselves. 

\begin{table}[h]
\begin{center}
\begin{tabular}{|l|r|r|r|}
\hline Smoothing strategy & W-cycle & V-cycle & $\omega$ \\ 
\hline \hline Unoptimized & 0.64 & 0.65 & 1.0 \\ 
\hline Exact block-wise solving & 0.11 & 0.18 & - \\ 
 \hline \hline Adaptive smoothing steps & 0.51 & 0.53 & 1.0 \\ 
 + Additional interface smoothing & 0.34 & 0.34 & 1.0 \\ 
\hline \hline Full damping & 0.59 & 0.56 & 1.15 \\ 
\hline Interior damping & 0.44 & 0.49 & 1.55 / 1.0 \\ 
 + Additional interface smoothing & 0.30 &  0.30 & 1.55 / 1.0 \\
\hline Adaptive interior damping & 0.46  & 0.51 & variable \\ 
  + Additional interface smoothing & 0.31 & 0.34 & variable \\
\hline \hline Combined Methods & 0.15 & 0.19 & 1.55 / 1.0 \\ 
\hline 
\end{tabular}
\end{center}
\caption{Measured convergence factors for different smoothing strategies} \label{convHalfSphereResults}
\end{table}

\subsubsection{Locally adaptive number of smoothing steps}
\label{Locallyadaptivenumberofsmoothingsteps}
LFA predictions can help to adjust the number of relaxation steps according to 
the element. Since the communication overhead should remain the same, we adjust the 
smoothing procedure locally. In order to individually increase the amount of smoothing work for a structured 
region, we add 4-color smoother iterations during smoothing steps. During execution 
of those iterations, the ghost points of the element are not updated. Thus, it is like smoothing 
a smaller Dirichlet boundary value problem. We decrease the smoothing work for a structured region 
by skipping smoothing during a (global) smoothing step. 

\begin{figure}
\centering
\caption{\label{fig:adaptSmoothing}Histogram of the number of 4-color iterations applied in the structured regions in the adaptive smoothing step version.}
{\includegraphics[width=0.8\textwidth]{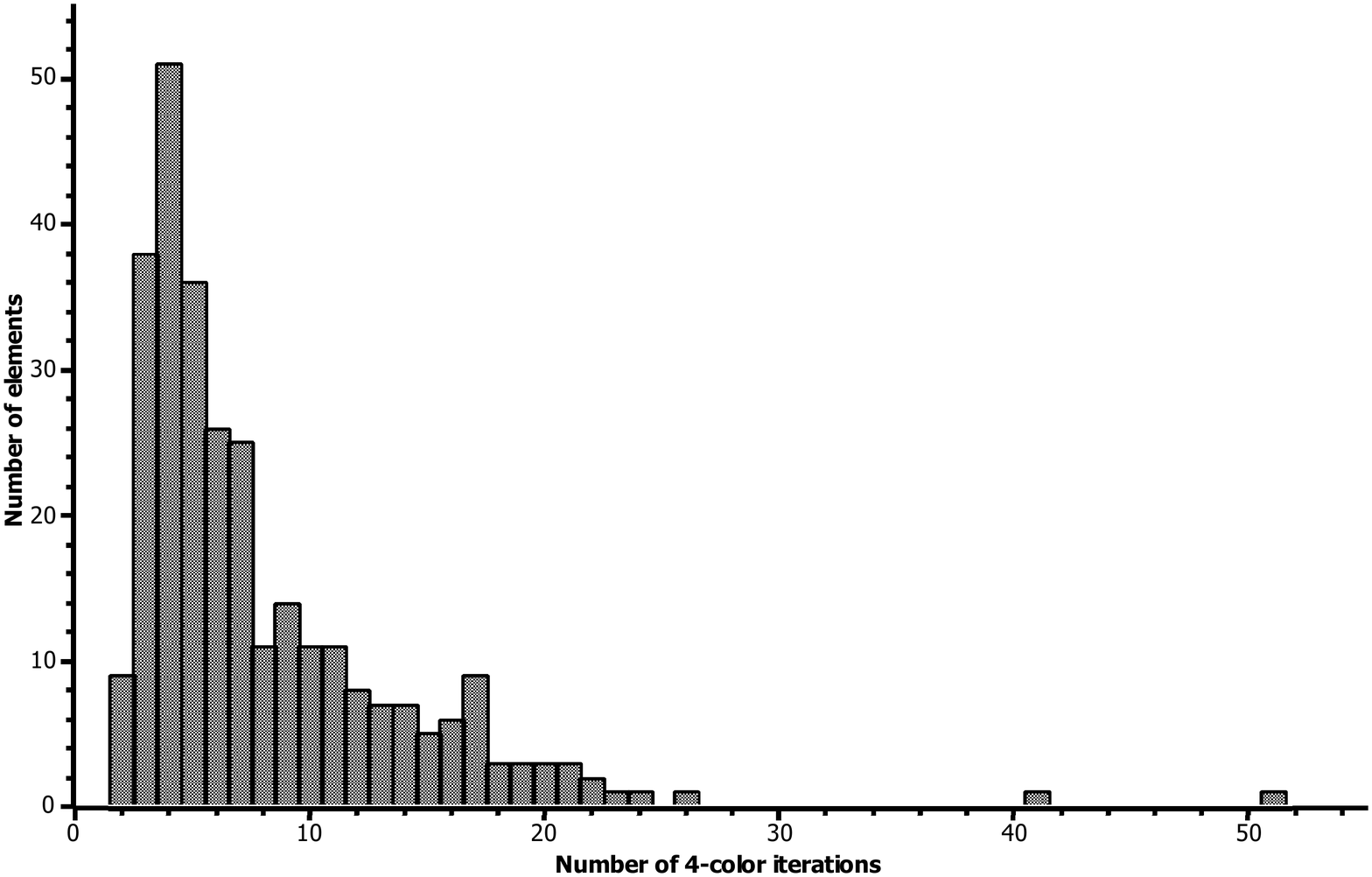}}\hfill
\end{figure}

For our experiment, the overall smoothing work should stay constant, thus it is more a 
redistribution of smoothing work. We start by setting the number of smoothing iterations 
for all structured regions to one. At this stage there is one real update of the ghost layers, 
but seven (three at pre- and four at post-smoothing) which would not be needed at all. 
Then all convergence rates are predicted by LFA and the remaining smoothing iterations are successively
distributed to the element with the highest convergence rate. After each step
the convergence rate of the updated element has to be recalculated. The smoothing iterations 
are distributed as equally as possible amongst between the updating steps of the ghost boundaries. The LFA evaluation
took around 5 seconds for our mesh in the setup time. But the advantage easily outweighs the required time. In 
table~\ref{convHalfSphereResults} \textit{adaptive smoothing steps} shows the effect of this optimization.
For our example, Figure~\ref{fig:adaptSmoothing} gives an impression of the smoothing work distribution.
For most elements $\nu_1+\nu_2 = 4$ steps are predicted to be sufficient, in an extreme cases even $\nu_1+\nu_2 = 51$ are predicted 
for an upper bound of $0.18$ on the convergence factor for each single structured region.

Moreover, we can add (local) smoothing iterations for the interface points, too. Since we do 
not save this additional work anywhere else, it is some overhead. However, it is neglectable compared 
to its possible effects (see \textit{+ additional interface smoothing} in Table~\ref{convHalfSphereResults}), 
since there are many more inner points than interface points. 

% It can be reasonable to set an upper limit for the number of smoothing steps. In our case we took $\nu = (25,25)$ for the LFA based prediction.

%With the LFA based estimation we can either reduce the number of
%V-cycles to 12, or perform 20 V-cycles with $\nu = (2,2)$. Figure
%\Fig{fig:convrates} plots the convergence rates of all three methods against
%the number of V-cycles.

%Using one core of an 2.8\,GHz \emph{Intel Core2 Quad} CPU, a V-cycle takes 16
%seconds. The setup phase of the LFA based approach takes 6 seconds. Thus,
%both optimization methods yield a large reduction of the time to solution.

\subsubsection{Locally adaptive damping}
\label{Locallyadaptivedamping}
In a first step we want to optimize the damping parameter. A simple choice 
is a constant damping parameter $\omega_{const} = \omega_1 = \omega_2 = 
\omega_3 = \omega_4$ for all points (see \textit{full damping} in Table~\ref{convHalfSphereResults}). By experiment 
$\omega_{const}=1.15$ showed up to be a good choice. Damping parameters around 
$1.25$ and more lead to divergence.  

%   \begin{figure}[!htb]
%   \caption{LFA estimated convergence factors for all elements without damping (origin of vectors), and its changing w.r.t. a damping parameter: }
%   \centering
%   \subfigure[single global damping parameter of $1.55$ and]{ \includegraphics[angle=-90, width=0.8\textwidth]{figures/constantDamping.pdf}} 
%   \subfigure[the arithmetic means of four optimal damping parameters per region.]{\includegraphics[angle=-90, width=0.8\textwidth]{figures/individualDamping.pdf}}
%   \label{fig:dampingOpt}
%   \end{figure}

%The LFA predicts optimal damping parameters between $1.30$ and $1.63$ (averaged over the four colors)
%for the different structured regions.
These damping parameters are low compared to LFA predicted optimal damping values in the interior of the
structured regions. 
Indeed, the interface points cause a problem in this case. Without damping the interface points we observe improved
convergence using a larger damping factor of $1.55$ (see \textit{interior damping}, Table~\ref{convHalfSphereResults}).
Increasing the experimentally found damping parameters on the macro vertices, edges and faces to $1.1$ nearly does not improve the convergence rate.
Similarly to the optimization in subsection~\ref{Locallyadaptivenumberofsmoothingsteps}, \textit{additional interface smoothing} 
is necessary to obtain the full potential of this optimization.
%Figure~\ref{fig:dampingOpt}~(a) shows the convergence rates before and after a 
%damping for each region of $1.55$. 

%Please note the y-axis has a logarithmic scale. The arrows in Figure~\ref{fig:dampingOpt}~(a)
%are sorted in the same way from left to right 
%like in Figure~\ref{fig:dampingOpt}~(b) for an easier comparison. Each arrow belongs to one structured region. Since all 
%convergence rates of the damped version are by far smaller than $1$, the interior of the structured regions might
%not be the reason for the divergence observed when using larger damping parameters. 

%For the interfaces between the structured regions it was
%important to choose lower damping parameter (in our case no damping) to prevent
%a diverging solution. Otherwise, we have to choose much smaller over-relaxation factors 
%and thus cannot obtain similar low convergence rates. 

%Increasing the damping parameters on the macro
%vertices, edges and faces to $1.1$ nearly does not improve the convergence rate.
%By experiment, a optimal damping parameter of $\omega = 1.65$ saves nine V-cycles.
%The asymptotic convergence factor improved to $0.49$.

To further improve the damping parameters, we tried to find optimal
relaxation parameters for each structured region. The four-color
damping parameters are determined by applying the downhill simplex method
\cite{Nel-1965} to LFA individually for each of the 293 elements on the coarsest mesh.
This is computationally quite expensive. The optimal 
damping parameter can vary quite strongly amongst the different stencils. Inside the
domain, quite low damping parameters e.g. $\omega = (1.43, 1.26, 1.17, 1.33)$
and large damping parameters e.g. $\omega = (1.75, 1.47, 1.46, 1.83)$ or
$\omega = (1.51, 2.00, 1.44, 1.46)$ are applied. However it did not 
improve the time to solution (see \textit{adaptive interior damping} in Table~\ref{convHalfSphereResults}). 
Exactly the same number of V-cycles like using 
the global damping parameter are required. The convergence factors are even a bit worse than using 
\textit{interior damping}. LFA predicts for some structured regions 
convergence factors down to $0.03$, while worst convergence factors are up to $0.51$ and $0.36$ 
with and without damping, respectively. 

Summarizing, the interface points benefit from additional smoothing. Further, the interior points 
can be damped more extremely, as long as damping is avoided on the interface points. It was not 
favorably to choose different damping parameter for the interior points.

\subsubsection{Combined methods}
In this section, we combine the various methods. As we will see, the positive 
effects w.r.t. the convergence factor will be additive to some extent. 
We chose the best variant of each previous subsection~\ref{Locallyadaptivedamping} and 
\ref{Locallyadaptivenumberofsmoothingsteps}.
Thus, the \textit{combined methods} of Table~\ref{convHalfSphereResults} include \textit{adaptive 
smoothing steps, interior damping}, and \textit{additional interface smoothing}.
%Two different sequences were tried out: Applying \textit{adaptive smoothing steps} first 
%and apply the damping parameters afterward, or choose the number of smoothing 
%steps by evaluating LFA to the damped smoothing. The second approach yields a better convergence. 
For this example the convergence rate is even similar to an \textit{exact 
block-wise solving}.

%% file: conclusion.tex
We were able to show a good correspondence between the LFA and HHG for
Jacobi, lexicographic and four-color Gauss-Seidel, line- as well as for plane-
type smoothers. The variation between the measured results of HHG and the LFA
predictions for the asymptotic convergence rate is below 10\% in all measured
cases. We presented a two-grid analysis for the Laplace operator discretized
by linear finite elements on tetrahedral grids. The four-color smoother requires
a decomposition of the Fourier space into minimal invariant subspaces, four and sixteen dimensional
for the smoothing  and  the two-grid analysis, respectively. 
The four-color smoother was the best choice for regular-like tetrahedra.
An optimal choice of individual damping parameters for the different colors has a stronger
impact to improve the convergence factor than it is the case for a lexicographic Gauss-Seidel.
Line-smoothing shows an especially good convergence for some tetrahedral
shapes, \eg{} for the wedge type. Plane-smoothing is very efficient for
needle or spindle shaped elements. For both multi-relaxation schemes, a zebra-wise updating
proved more effective than a lexicographic updating.

The agreement between the LFA and HHG can be seen as a validation for both
implementations. Furthermore, benefits of more complex solvers can first be
tested with the LFA before implementing them in HHG. This prevents development
effort that would be futile, if the solver turned out not to be suitable for
HHG.

Considering semi-structured grids, the interfaces between the structured blocks are often used, \eg{}, as boundaries
for a parallelization. Thus points on the boundaries have to be updated before or
after the structured parts to keep the communication low. While for point-wise
smoothing this proceeding has nearly no influence on the convergence, this is
not the case for line- or plane-smoothing. We experienced a strong effect when
interfaces separated a collective update in several parts. While an additional
smoothing of the interface points and over-relaxation can slightly reduce this effect,
further work has to be done for an improvement. On the other hand we could show, with
three different smoothers on three degenerated tetrahedral regions, that there is no
inherent problem for locally adaptive smoothing. It can be very advantageous to use
different smoothing strategies for the different sub-domains.

Because of the refinement strategy in HHG, child elements of a tetrahedron have
the same anisotropy as the parent element. Thus, it can be worthwhile doing some
calculations in order to find out good smoothing parameters, since these
parameters are applied to a huge number of tetrahedrons on the finer grids. In
the last part of the paper we focused on the number of pre-/post-smoothing steps 
and local adaptively chosen damping parameters. It turned
out that by redistributing the number of smoothing steps and applying different 
local damping parameters over a semi-structured
test domain can significantly reduce the convergence.

% An open question is how these improvments behave in a parallel execution,
% since they are currently only working locally on each compute node. We expect
% that the effects of the smoothing step redistribution decreases with increasing
% number of processors. However, the presented methods have a good potential for
% load balancing since they are based on a~priori estimations.

%% file: biblio.tex
%\biboptions{longnamesfirst,angle,semicolon}
\bibliographystyle{model1-num-names}